\documentclass[11pt,reqno]{amsart}
\usepackage{amsmath,amstext,amsfonts,amsthm,amssymb,dsfont}
\usepackage[all]{xy}
\usepackage{eucal}
\usepackage{graphicx}
\swapnumbers
\newtheorem{thm}[subsubsection]{Theorem}
\newtheorem{lem}[subsubsection]{Lemma}

{  \theoremstyle{definition}

}
\numberwithin{equation}{subsection}
\newcommand{\iso}{\overset{\sim}{\longrightarrow}}

\newcommand{\shriekotimes}{\otimes^{!}}

\newcommand{\bea}{\begin{eqnarray*}}
\newcommand{\eea}{\end{eqnarray*}}
\newcommand{\bean}{\begin{eqnarray}}
\newcommand{\eean}{\end{eqnarray}}

\newcommand{\fg}{\mathfrak g}

\newcommand{\CA}{\mathcal{A}}

\newcommand{\CC}{\mathcal{C}}
\newcommand{\CD}{\mathcal{D}}

\newcommand{\CL}{\mathcal{L}}
\newcommand{\CM}{\mathcal{M}}
\newcommand{\CN}{\mathcal{N}}
\newcommand{\CO}{\mathcal{O}}
\newcommand{\CP}{\mathcal{P}}

\newcommand{\CS}{\mathcal{S}}
\newcommand{\CT}{\mathcal{T}}
\newcommand{\CV}{\mathcal{V}}

\newcommand{\BC}{\mathbb{C}}

\newcommand{\BZ}{\mathbb{Z}}
\newcommand{\Vac}{\mathds{1}}

\newcommand{\nc}{\newcommand}

\nc{\Id}{\text{Id}}
\nc{\la}{\lambda}

\begin{document}


\newcommand{\Real}{\mathbb R}
\newcommand{\HH}{\mathbb H}
\newcommand{\QQ}{\mathbb Q}
\newcommand{\ZZ}{\mathbb Z}
\newcommand{\LL}{\mathbb L}
\newcommand{\VV}{\mathbb V}
\newcommand{\MM}{\mathbb M}
\newcommand{\PP}{\mathbb P}
\newcommand{\RR}{\mathbb R}

\newcommand{\cA}{\mathcal{A}}
\newcommand{\cB}{\mathcal{B}}
\newcommand{\cC}{\mathcal{C}}
\newcommand{\cext}{\cC ext}
\newcommand{\cD}{{\mathcal{D}}}
\newcommand{\cE}{{\mathcal{E}}}

\newcommand{\cF}{{\mathcal{F}}}
\newcommand{\cG}{{\mathcal{G}}}
\newcommand{\cH}{{\mathcal{H}}}
\newcommand{\cJ}{{\mathcal{J}}}

\newcommand{\cL}{{\mathcal{L}}}
\newcommand{\cM}{{\mathcal{M}}}
\newcommand{\cN}{\mathcal{N}}
\newcommand{\cO}{{\mathcal{O}}}
\newcommand{\cP}{{\mathcal{P}}}
\newcommand{\cQ}{{\mathcal{Q}}}
\newcommand{\cR}{{\mathcal{R}}}
\newcommand{\cS}{{\mathcal{S}}}
\newcommand{\cT}{{\mathcal{T}}}
\newcommand{\cV}{{\mathcal{V}}}
\newcommand{\cW}{{\mathcal{W}}}
\newcommand{\vext}{\cV ext}

\newcommand{\cZ}{{\mathcal{Z}}}

\title{strongly homotopy chiral  algebroids}

\author{ F.Malikov }\thanks{partially supported by an NSF grant}
\maketitle

\begin{abstract}

We introduce and classify the objects that appear in the title of the paper.

 \end{abstract}

\section{introduction}
\label{intro}
\subsection{ } An old paper by B.~Feigin and A.~Semikhatov, \cite{FS}, suggests the following construction and proves the following
theorem, a result of rare beauty. Start with the Koszul resolution $K(\BC[x],\{x^m\})$
\[
0\longrightarrow\BC[x]\stackrel{x^m}{\longrightarrow}\BC[x]\longrightarrow 0.
\]
Now, ``chiralize'' this resolution. Namely, consider the vertex algebra $K^{ch}(\BC[x],\{x^m\})$ that is 
generated by two pairs of fields, $\partial_x(z)$, $x(z)$, both even,
and $\partial_\xi(z)$, $\xi(z)$, both odd, such that the only nontrivial commutation relations are as follows
\[
[\partial_x(z),x(w)]=[\partial_\xi(z),\xi(w)]=\delta(z-w).
\]
Give this vertex algebra the differential 
\[
D=\oint x(z)^m\partial_\xi(z)\,dz.
\]
The aforementioned result asserts that the cohomology $H_D(K^{ch}(\BC[x],\{x^N\}))$ is a direct sum of $m$ distinct unitary
representations of the celebrated $N=2$ superconformal Lie algebra generated by the classes of the classical Koszul cohomology
$1, x, x^2,\ldots, x^{m-1}$.

Some later work,\cite{BD,MSV,GMS}, makes it clear that this cohomology  has the meaning of of an algebra of chiral differential
operators over a fat point, $\text{Spec}\BC[x]/(x^m)$, a highly singular affine scheme.

We would like to understand whether the vertex algebra $H_D(K^{ch}(\BC[x],\{x^N\}))$ can be defined conceptually, and not by writing
formulas.

\subsection{ } These notes are then about the following circle of ideas. Let $A$ be a smooth affine $\BC$-algebra, $T_A=\text{Der}(A)$ the corresponding tangent Lie algebroid. A Picard-Lie $A$-algebroid is
an exact sequence
\[
0\longrightarrow A\stackrel{\iota}{\longrightarrow} \CL\stackrel{\sigma}{\longrightarrow} T_A\longrightarrow 0,
\]
where $\CL$ is a Lie $A$-algebroid, and the arrows respect all the structures. The category of Picard-Lie algebroids is governed by the truncated De Rham complex
$\Omega_A^1\rightarrow\Omega_A^{2,cl}$, \cite{BB}. Informally speaking, deformations of the bracket involve closed 2-forms, and so objects are labeled by closed 2-forms, 
morphisms by 1-forms; formally, the category of Picard-Lie
algebroids is a $\Omega_A^{[1,2>}$-torsor, where $\Omega_A^{[1,2>}$ is a category with objects $\Omega^{2,cl}_A$ and morphisms defined by $Hom(\beta,\gamma)=
\{\alpha\in \Omega^1_A\text{ s.t. }d\alpha=\beta-\gamma\}$.

A chiral $A$-algebroid is an exact sequence
\[
0\longrightarrow J_\infty A\stackrel{\iota}{\longrightarrow} \CL^{ch}\stackrel{\sigma}{\longrightarrow} J_\infty T_A\longrightarrow 0,
\]
where 

$J_\infty A$ is the corresponding jet-algebra, and in particular, a commutative chiral algebra; 

$J_\infty T_A$ is a tangent Lie* algebroid, which is an analogue of $T_A$
in the world of the Beilinson-Drinfeld pseudo-tensor categories \cite{BD}; in particular, $J_\infty T_A$ carries the compatible structures of a  $J_\infty A$-module
and a Lie* algebra;

$\CL^{ch}$ is a Lie* algebra and a {\em chiral } $J_\infty A$-module (this last notion is different from that of an ordinary $J_\infty A$-module used a line above, and this has consequences);

the morphisms $\iota$ and $\sigma$ respect all the structures.

 It is appropriate at this point to make a terminological remark. This paper could not have been written
outside the framework created in \cite{BD}; however, our emphasis (for simplicity and as  a reflection of personal limitations) is entirely on translation-invariant objects over
$\BC$, and so our chiral algebras typically are vertex algebras, Lie* algebra are vertex Lie algebras, etc. The objects we are dealing with, instead of being $D$-modules, are
$R$-modules, $R=\BC[\partial]$, $\partial$ being thought of as  $d/dx$. There is little doubt that most of our discussion apply in greater generality.

Be it as it may, the category of chiral $A$-algebroids is a torsor over $C^{[1,2>}_{J_\infty A}(J_\infty T_A,J_\infty A)$ \cite{BD}, where $C^{\bullet}_{J_\infty A}(J_\infty T_A,J_\infty A)$
is a De Rham-Chevalley complex,  an object introduced in \cite{BD} and which is  to $J_\infty T_A$  what $\Omega^\bullet_A$ is to $T_A$.  
A much less general but more explicit result was proved in \cite{GMS}.

An attempt to deal with a singular $A$ leads to an $A$ that is  a finitely generated, polynomial DG algebra, with differential of degree 1. In both of the cases just considered the exact
sequences still make perfect sense in the category of the corresponding DG objects, and so do the complexes, such as $\Omega_A^\bullet$ or 
$C^{\bullet}_{J_\infty A}(J_\infty T_A,J_\infty A)$, which now acquire an extra grading and  differential. However, classification of such exact sequences again leads
to the familiar truncated complexes, such as $\Omega_A^1\rightarrow\Omega_A^{2,cl}$, which homotopically makes little sense.

It appears that the right thing to do is to define a Picard-Lie $\infty$-algebroid, where $\CL$ fits the above exact sequence and is allowed to be an $A$-module and
a $\text{Lie}_\infty$ algebra, but not necessarily an ordinary Lie algebra. We prove (Lemma~\ref{res-class-picardlie-infty}) that the category of such algebroids is a torsor over $L\Omega^{[1,2>}_A$, an analogue
of $\Omega^{[1,2>}_A$ where the usual De Rham complex is replaced  with the derived De Rham introduced by Illusie \cite{Ill1,Ill2}. Informally, the fact that a $\text{Lie}_\infty$
algebra carries an infinite family of higher brackets creates an avenue for deformations by higher total degree   2 forms $\alpha_n\in\Omega^{n}_A[2-n]$. Similarly, the fact that a morphism
of a $\text{Lie}_\infty$ algebra is a morphism of the corresponding symmetric algebra provides for morphisms determined by higher total degree 1 forms
$\alpha_n\in\Omega^{n}_A[1-n]$.

\begin{sloppypar}
Further, we introduce the concepts of a $\text{Lie}^*_\infty$ algebra (sect.~\ref{lie-star-infty}) and
 of a chiral $\infty$-algebroid (sect.~\ref{descr-obf-chiral-infty-al}) by allowing $\CL^{ch}$ in the above exact sequence to be a chiral $J_\infty A$-module and
a $\text{Lie}^*_\infty$ algebra.
The main result of this paper consists in the  classification of chiral $\infty$-algebroids, Theorem~\ref{class-oo-chi-alge}. These form a torsor over
the Illusie-De Rham-Chevalley $LC^{[1,2>}_{J_\infty A}(J_\infty T_A,J_\infty A)$.  The appearance of higher total degree 2 forms is easy to anticipate, but does require the use of the
Beilinson-Drinfeld *-operations. An interpretation of a $\text{Lie}^*_\infty$ algebra as a symmetric algebra, which is a requirement for higher morphisms, is perhaps the
least familiar part of the present work and uses the Beilinson-Drinfeld category $Mod(R^{\CS})$,  a ``tensor enveloping category'' of the pseudo-tensor category
of *-operations \cite{BD}; this is done in sects.~\ref{from pseudo-tensor-to tensor}, \ref{liestarinfty-via-coderiv}.
\end{sloppypar}

\subsection{ } The main results  are stated and proved in sects.~\ref{picard-lie-infty-alg} and \ref{chiral-infty-algebro}; it is for the sake of these
sections that the paper was written and might be read.  The purpose of the rest is  to facilitate the references.   This is especially true of
sect.~\ref{beilinson-drinfeld}, which  can be characterized as directed at a ``VOA insider untrammelled by algebro-geometric affections''
\cite{BD}. Sect.~\ref{cdo} contains a reminder on algebras of chiral differential operators in the generality used in papers
such as \cite{MSV,GMS}, but in the form suggested  by \cite{BD}; the clarity achieved using the latter approach is quite striking.

\subsection{ } These notes originated in an attempt to understand the mysterious unpublished manuscript by V.Hnich, \cite{Hin}. It would be fair
if V.Hinich were an author, but he refused. I am grateful to V.Hinich for sharing his ideas with me.

\section{TDO}
\label{TDO}
\subsection{ }
\label{defn-of -tdo}
Let $A$ be a commutative unital  $\BC$-algebra, $T_A$ the Lie algebra of derivations of $A$. The graded symmetric algebra $S^\bullet_A T_A$ is naturally
a Poisson algebra. An algebra $\CD^{tw}_A$ is called an algebra of twisted differential operators over $A$, TDO for brevity, if it carries a filtration
$F_0(\CD^{tw}_A)=A\subset\cdots\subset F_n(\CD^{tw}_A)\subset\cdots$, $\bigcup_n 
F_n(\CD^{tw}_A)=\CD^{tw}_A$, s.t. the corresponding
graded object is isomorphic to $S^\bullet_A T_A$ is a Poisson algebra. 

In a word, a TDO is a quantization of $S^\bullet_A T_A$.

\subsection{ }
\label{defn-defn-p-lie}
The key to classification of TDOs is the concept of a Picard-Lie $A$-algebroid.  $\CL$ is called a Lie $A$-algebroid if it is
 a Lie algebra,  an $A$-module, and is equipped with anchor, i.e., a Lie algebra and an $A$-module map
 $\sigma: \CL\rightarrow T_A$ s.t. the $A$-module structure map
 \begin{equation}
 \label{lie-alg-comp-str}
 A\otimes \CL\longrightarrow\CL
 \end{equation}
 is an $\CL$-module morphisms. Explicitly,
 \begin{equation}
 \label{brac-vect-fields}
 [\xi,a\tau]=\sigma(\xi)(a)\tau+a[\xi,\tau],\; a\in A,\, \xi,\tau\in \CL.
 \end{equation}
 A Picard-Lie $A$-algebroid is a Lie $A$-algebroid $\CL$ s.t. the anchor fits in an exact sequence
\begin{equation}
\label{defn-pl-alg-exact}
0\longrightarrow A\stackrel{\iota}{\longrightarrow} \CL\stackrel{\sigma}{\longrightarrow} T_A\longrightarrow 0,
\end{equation}
where the arrows respect all the structures involved; in particular, $A$ is regarded as an $A$-module and an abelian Lie algebra, and $\iota$
makes it an $A$-submodule and an abelian Lie ideal of $\CL$. Furthermore, the induced action of  $T_A=\CL/A$ on $A$ must be equal to the canonical action of 
$T_A$ on $A$.

Morphisms of Picard-Lie $A$-algebroids are defined in an obvious way to be morphisms of exact sequences (\ref{defn-pl-alg-exact}) that preserve all  the structure
involved. Each such morphism is automatically an isomorphism and we obtain a groupoid $\CP\CL_A$.

\subsection{ }
\label{defn-p-lie}Classification of Picard-Lie $A$-algebroids that split as $A$-modules is as follows. We have a canonical such algebroid, $A\oplus T_A$ with bracket
\[
[a+\xi,b+\tau]=\xi(b)-\tau(a)+[\xi,\tau].
\]
Any other bracket must have the form
\[
[\xi,\tau]_{new}=[\xi,\tau]+\beta(\xi,\tau),\; \beta(\xi,\tau)\in A.
\]
The $A$-module structure axioms imply that $\beta(.,.)$ is $A$-bilinear, the Lie algebra axioms imply that, in fact,
$\beta\in\Omega^{2,cl}_A$. Denote this Picard-Lie algebroid by $T_A(\beta)$. Clearly, any Picard-Lie $A$-algebroid is isomorphic to
$T_A(\beta)$ for some $\beta$.

A morphism $T_A(\beta)\rightarrow T_A(\gamma)$ must have the form $\xi\rightarrow\xi+\alpha(\xi)$ for some $\alpha\in\Omega^1_A$.
A quick computation will show that
\[
Hom(T_A(\beta),T_A(\gamma))=\{\alpha\in\Omega^1_A\text{ s.t. } d\alpha=\beta-\gamma\}.
\]
This can be rephrased as follows. Let $\Omega^{[1,2>}_A$ be a category with objects $\beta\in\Omega^{2,cl}_A$, morphisms
$Hom(\beta,\gamma)=\{\alpha\in\Omega^1_A\text{ s.t. } d\alpha=\beta-\gamma\}$. Then the assignment $(\gamma,T_A(\beta)\mapsto T_A(\beta+\gamma)$
defines a categorical action of  $\Omega^{[1,2>}_A$ on $\CP\CL_A$ which makes $\CP\CL_A$ into  an $\Omega^{[1,2>}_A$-torsor. The isomorphism classes
of this catetgory are in 1-1 correspondence with the De Rham cohomology $\Omega^{2,cl}_A/d\Omega^1_A$, and the automorphism group of an object is $\Omega^{1,cl}_A$.

\subsection{ }
\label{pic-li-varie}
If $X$ is a smooth algebraic variety, then the above considerations give the category of Picard-Lie algebroids over $X$,
$\CP\CL_X$, which is a torsor over $\Omega^{[1,2>}_X$ or, perhaps, a gerbe bound by the sheaf complex
$\Omega^1_X\rightarrow \Omega^{2,cl}_X$. This gerbe has a global section, the standard $\CO_X\oplus\CT_X$. The isomorphism of classes of such algebroids are
in 1-1 correspondence with the cohomology group $H^1(X,\Omega^1_X,\rightarrow\Omega^{2,cl}_X)$ ($\Omega^1_X$ being placed in degree 0), and the automorphism group of an object is $H^0(X,\Omega^{1,cl}_X)$.

\subsection{ }
\label{lie-alg-envelope}
The concept of the universal enveloping algebra of a Lie algebra has a Lie algebroid version, which reflects a partially defined multiplicative structure on $\CL$.

Let $F(\CL)$ be a free unital associative $\BC$-algebra generated be the Picard-Lie $A$-algebroid $\CL$ regarded as a vector space over $\BC$. We denote
by $\ast$ its multiplication and by ${\bf 1}$ its unit. Define the  universal enveloping algebra $U_A(\CL)$ to be the quotient of $F(\CL)$ be the ideal generated by
the elements $\xi\ast\tau-\tau\ast\xi-[\xi,\tau]$,   $a\ast\xi -a\xi$, ${\bf 1}-1_A$, where $1_A$ is the unit of $A$.

It is rather clear that  $U_A(\CL)$ is a TDO (sect.~\ref{defn-of -tdo}), and the assignemnt $\CL\mapsto U_A(\CL)$ is an equivalence of categories if $A$ is smooth, i.e., if
$MaxSpec(A)$ is a smooth affine variety.

\section{Picard-Lie $\infty$-algebroids.}
\label{picard-lie-infty-alg}
\subsection{ }
\label {lie-infty}
Let $L$ be a graded vector space s.t. $L=\oplus_n L_n$. Call a map $f: L^{\otimes n}\rightarrow L$ {\em antisymmetric} if 
$f(x_{\sigma_1},\ldots,f_{\sigma_n})=sgn(\sigma)\epsilon(\sigma,\vec{x})f(x_1,\ldots,x_n)$, where 
\[
\epsilon(\sigma,\vec{x})=\prod_{i<j\text{ s.t. }\sigma_i>\sigma_j}(-1)^{x_{\sigma_i}x_{\sigma_j}}.
\]
Similarly, a map $f: L^{\otimes n}\rightarrow L$ is called {\em symmetric} if $f(x_{\sigma_1},\ldots,f_{\sigma_n})=\epsilon(\sigma,\vec{x})f(x_1,\ldots,x_n)$.

The space $L^{\otimes n}$ is naturally graded and we say that $f: L^{\otimes n}\rightarrow L$ has degree $m$ if  $f((L^{\otimes n})_n)\subset L_{n+m}$.

A $\text{Lie}_\infty$ algebra (cf. \cite{LM}) is a graded vector space $L$ with a collection of antisymmetric maps $l_n: L^{\otimes n}\rightarrow L$, $n=1,2,3,...$, s.t. $\text{deg}l_n=2-n$
and the following identity is satisfied for each $k\geq 1$
\begin{equation}
\label{jacobi-for-infty-ordinary}
\sum_{i+j=k+1}\sum_{\sigma}\text{sgn}(\sigma)\epsilon(\sigma,\vec{x})(-1)^{i(j-1)}l_j(l_i(x_{\sigma_1}\ldots x_{\sigma_i}),x_{\sigma_{i+1}}\ldots x_{\sigma_n})=0,
\end{equation}
where $\sigma$ runs through the set of all $(i,n-i)$ unshuffles, i.e.,  $\sigma\in S_n$ s.t. 
$\sigma_1<\sigma_2<\cdots<\sigma_i$ and $\sigma_{i+1}<\sigma_{i+2}<\cdots<\sigma_n$.

A {\em strict } $\text{Lie}_\infty$ algebra morphism from $(L,\{l_n\})$ to $(L',\{l_n'\})$ is a degree 0 map $f: L\rightarrow L'$ s.t.
 $f\circ l_n=l_n'\circ f^{\otimes n}$ for each $n\geq 1$.
 
 \subsection{ }
 \label {def-via-coderiv}
 Let $S^\bullet L$ be the free (graded) commutative algebra generated by $L$, i.e., the quotient of  the tensor algebra $T(L)=\oplus_n L^{\otimes n}$ by the 2-sided
 ideal generated by the elements $x\otimes y-(-1)^{xy}y\otimes x$. In what follows the class of $x\otimes y\otimes z\cdots$ in $S^\bullet L$ is denoted by $xyz\cdots$.
 
$S^\bullet L$ carries a coalgebra structure $\Delta: S^\bullet L\rightarrow S^\bullet L\otimes S^\bullet L$ defined by
\begin{equation}
\label{def-coalg-usu-sy-alg}
\Delta(x_1x_2\cdots x_n)=\sum_i\sum_\sigma \epsilon(\sigma,\vec{x})x_{\sigma_1}x_{\sigma_2}\cdots x_{\sigma_i}\otimes x_{\sigma_{i+1}}\cdots x_{\sigma_n},
\end{equation}
where the summation is extended to all $(i,n-i)$-unshuffles $\sigma$.

A {\em coderivation} is a linear map $f: S^\bullet L\rightarrow S^\bullet L$ s.t. $\Delta\circ f=(f\otimes 1+1\otimes f)\circ\Delta$, where
 $(1\otimes f)(x\otimes y)=(-1)^{xf}x\otimes f(y)$ -- the Koszul rule. The space of all coderivations is a Lie subalgebra of $\text{End}(S^\bullet L)$.

The Lie algebra of all coderivations of $S^\bullet L$ that preserve the filtration by degree is identified with $Hom_\BC(S^\bullet L,L)$ where $f\in Hom_\BC(S^nL,L)$ defines the coderivation
\begin{equation}
\label{form-fo-code-ordina-ca}
f(x_1\cdots x_N)=\sum_\sigma \epsilon(\sigma,\vec{x})f(x_{\sigma_1},\cdots x_{\sigma_n})x_{\sigma_{n+1}},...,x_{\sigma_N},
\end{equation}
the summation being extended to all $(n,N-n)$-unshuffles $\sigma$.

Along with $L$, consider $L[1]$, the graded space s.t. $L[1]_n=L_{n+1}$.  Denote by $s$ the identity map $L[1]\rightarrow L$; its degree is 1. There arises a map
\[
s^{\otimes n}: L[1]^{\otimes n}\longrightarrow L^{\otimes n},\; x_1\otimes x_2\otimes\cdots x_n\mapsto
(-1)^{(n-1)x_1+(n-2)x_2+\cdots +x_{n-1}}x_1\otimes x_2\otimes\cdots x_n,
\]
where the sign is forced upon us by the Koszul rule

Given $f\in Hom_\BC(L^{\otimes n},L)$ define $\hat{f}=(-1)^{n(n-1)/2}s^{-1}\circ f\circ s^{\otimes n}\in Hom_\BC(L[1]^{\otimes n},L[1])$. If $f$ is antisymmetric, then $\hat{f}$ is symmetric,
hence defines an element of $Hom_\BC(S^n L[1], L[1])$. The latter map as well as the corresponding coderivation of $S^\bullet L[1]$ will also be denoted by $\hat{f}$.

Any $\text{Lie}_\infty$ algebra $(L,\{l_n\})$ defines, therefore, a coderivation, $\sum_{n\geq 1}\hat{f}_n$, of $S^\bullet L[1]$. A well-known result, \cite{LM,LS}, asserts
that 

{\em this construction sets up a 1-1 correspondence between $\text{Lie}_\infty$ algebra structures on $L$ and coderivations of  $S^\bullet L[1]$ of degree 1 and square 0.}

This result prompts the following definition, \cite{LM}.

Define a $\text{Lie}_\infty$ algebra {\em morphism} $f: (L,\{l_n\})\rightarrow (M,\{m_n\})$ to be a morphism of coalgebras
with derivations, i.e., a coalgebra morphism $f: S^\bullet L[1]\rightarrow S^\bullet M[1]$ s.t. $f\circ\sum_n\hat{l}_n=\sum_n\hat{m}_n\circ f$.

Notice that a coalgebra morphism $f: S^\bullet L[1]\rightarrow S^\bullet M[1]$ is a collection of degree 0 symmetric maps $f_n:S^n L[1]\rightarrow  M[1]$, $n\geq 1$,
s.t.
\begin{multline}
\label{expl-descr-morph-infty}
f(x_1\cdots x_n)=\\
\sum_{1\leq i_1<i_2<i_3<\cdots <n}\sum_\sigma \epsilon(\sigma,\vec{x})f_{i_1}(x_{\sigma_1},...,x_{\sigma_{i_1}})f_{i_2-i_1}(x_{\sigma_{i_1+1}},...,x_{\sigma_{i_2}})
\cdots f_{n-i_{k}}(x_{\sigma_{i_{k}+1}},...,x_{\sigma_{n}}),
\end{multline}
where the summation is extended over all those $\sigma$ that satisfy $\sigma_1<\sigma_{i_1+1}<\sigma_{i_2+1}<\cdots$ and $\sigma_{i_j+\alpha}<\sigma_{i_j+\beta}$
as long as $i_j+\alpha< i_j+\beta\leq i_{j+1}$

\subsubsection{Remark.}
\label{remark-on-tens-vs-assinfty-alg}
We would like to recall, for future reference, that similar and simpler formulas can be written for the tensor algebra, $TL$, in place of the symmetric algebra. It is also
a coalgebra with comultiplication
\begin{equation}
\label{defn-comul-us-te-alge}
\Delta: TL\longrightarrow TL\otimes TL,\; x_1 x_2\cdots x_n\mapsto\sum_{i=0}^{n}x_1 \cdots x_i\otimes x_{i+1}\cdots x_n.
\end{equation}
The Lie algebra of coderivations of $TL$ that preserve the filtration by degree is identified with $Hom_\BC(TL,L)$, where  $f\in Hom_{\BC}(T^mL,L)$ defines the
coderivation
\[
f(x_1\otimes\cdots \otimes x_n)=\sum_{i=1}^{n-m}(-1)^{f\cdot(x_1+\cdots +x_i)}x_1\otimes\cdots x_i\otimes f(x_{i+1}\otimes\cdots\otimes x_{i+m})\otimes\cdots\otimes x_n.
\]
Similarly to the Lie case, one can define the concept of an $Ass_\infty$-algebra and verify that this structure on $L$ is the same thing as  a degree 1 and square 0 coderivation
of $TL[1]$, \cite{LM}.

\subsection{ }
\label{infty-pl-algebr}

Let $(A,\partial)$ be a differential graded algebra, differential  $\partial$ having degree 1.   The tangent $A$-algebroid $T_A$ is then a a differential graded Lie algebra,
hence a $\text{Lie}_\infty$ algebra with $l_1(.)=[\partial,.]$, $l_2(.,.)=[.,.]$, $l_n=0$ if $n\geq 3$.

A Picard-Lie $\infty$-algebroid is an exact sequence, cf. sect.~\ref{defn-p-lie},
\begin{equation}
\label{defn-pic-lie-infty}
0\longrightarrow A\stackrel{\iota}{\longrightarrow}\CL_A\stackrel{\sigma}{\longrightarrow} T_A\longrightarrow 0,
\end{equation}
where 

(i) $\CL_A=(\CL_A,\{m_n\})$ is a (differential graded) $A$-module and a $\text{Lie}_\infty$ algebra, $\sigma$, the anchor map, is a morphism of $A$-modules and a {\em strict}
morphism (sect.~\ref{lie-infty}) of $\text{Lie}_\infty$ algebras; 

(ii) the $\text{Lie}_\infty$ algebra structure of $\CL_A$ satisfies
\begin{eqnarray}
\label{mod-alg-comp-infty}
m_1(ax)&=&\partial(a) x+(-1)^a a m_1(x),\nonumber\\
 m_2(y,ax)&=&[\sigma(y),a] x+(-1)^{ya}a m_2(y,x),\\ 
m_j&\text{is}&A\text{-linear  and factors through  $\sigma$  if }j\geq 3.\nonumber
\end{eqnarray}


(iii) the embedding $\iota: A\rightarrow\CL_A$ is a morphism of $A$-modules and a strict morphism of $\text{Lie}_\infty$ algebras, where $A$ is regarded as an abelian such algebra
($l_2=l_3=\cdots=0$), with $l_1=\partial$.

An example is provided by the ordinary Picard-Lie algebroid $A\oplus T_A$ with differential $[\partial,.]$ and the obvious $\iota$ and $\sigma$.

\subsection{ }
\label{class-picard-lie-infty-ord}

Define a morphism of Picard-Lie $\infty$-algebroids, $ (\CL_A,\iota,\sigma)\rightarrow(\CL_A',\iota',\sigma')$ to be a morphism of $\text{Lie}_\infty$ algebras
$f=\{f_n\}: (\CL_A,\{l_n\})\rightarrow (\CL_A',\{l_n'\})$ (as defined in sect.~\ref{def-via-coderiv}) that satisfies the following 2 conditions:

(i) each $f_n$, which is a map $S^n\CL_A[1]\rightarrow \CL_A'[1]$, is $A$ $n$-linear; furthermore, if $n>1$, then $f_n$ factors through a map
$S^n(\CL_A/A)[1]\stackrel{\sigma}{\longrightarrow}S^nT_A[1]\longrightarrow A[1]\hookrightarrow \CL_A'$;

(ii) the component $f_1:\CL_A\rightarrow\CL_A'$ (see (\ref{expl-descr-morph-infty})
makes the following a commutative diagram of $A$-module morphisms:

\begin{displaymath}
\label{comm-diagr}
\xymatrix{ 0\ar[r]&A\ar[r]^{\iota}&\CL_A\ar[r]^{\sigma}\ar[d]^{f_1}&
T_A\ar[r]&0\\
 0\ar[r]&A\ar@{=}[u] \ar[r]^{\iota'}& \CL_A'\ar[r]^{\sigma'}&T_A\ar@{=}[u]\ar[r]&0}
 \end{displaymath}

\subsection{ }
\label{inftya-alg-categ}
In order to describe the category of Picard-Lie $\infty$-algebroids, 
consider the  De Rham complex $\hat{\Omega}^\bullet_A[\bullet]$. It is the usual De Rham complex of $A$, which is now a bi-complex with extra differential $[\partial,.]=\text{Lie}_\partial$
and extra grading coming from that of $A$,
except that it is completed in the De Rham direction (see \cite{Ill1,Ill2}.) Denote by $L\Omega_A^\bullet$ the corresponding total complex s.t. 
$L\Omega_A^k=\prod_{n=1}^\infty\Omega^n_A[2-n]$ and the differential equals $d+\partial$.  

As in sect.~\ref{defn-p-lie}, define the category $L\Omega_A^{[1,2>}$ with objects elements of $L\Omega_A^2$ and morphisms defined by
$Hom(\beta,\gamma)=\{\alpha\in L\Omega^1_A\text{ s.t. } (d+\partial)\alpha=\beta-\gamma\}$.

\begin{lem}
\label{res-class-picardlie-infty}
The category of Picard-Lie $\infty$-algebroids over $A$ is a torsor over the category $L\Omega^{[1,2>}_A$, cf. sect.~\ref{defn-p-lie},\ref{pic-li-varie}.
\end{lem}

Indeed, having chosen the standard $A\oplus T_A$ as a reference point any other such algebroid is obtained by replacing
the standard $l_j$ with $\tilde{l}_j=l_j+\alpha_j$. The definition, sect.~\ref{infty-pl-algebr}, implies $\alpha_j\in\Omega_A^j[2-j]$. Among the terms of the relations (\ref{jacobi-for-infty-ordinary}) 
for $\{\tilde{l}_n\}$ the nontrivial ones will be of either of the forms $l_2\circ l_{j-1}$, $l_{j-1}\circ l_2$, $l_1\circ l_{j}$, $l_{j}\circ l_1$. The first two will give
$d_{DR}(\alpha_{j-1})$, the last two will give $\text{Lie}_\partial(\alpha_j)$. Overall, one obtains the cocycle condition: $(d_{DR}+\text{Lie}_\partial)(\{\alpha_n\})=0$.

Morphisms are  collections of degree 0 $A$-multilinear maps $\hat{\beta}_n: S^n T_A[1]\rightarrow T_A[1]$. The actual morphism they define operates as follows:
\[
f(x)=x+\hat{\beta}_1(x),
\]
\[
f(x_1,x_2)=(x_1+\hat{\beta}_1(x_1))(x_2+\hat{\beta}_1(x_2))+\hat{\beta}_2(x_1,x_2),
\]
etc., cf.(\ref{expl-descr-morph-infty}).  Such morphisms are automatically automorphisms; in fact, $f^{-1}$ is the morphism defined by the collection $\{-\hat{\beta}_n\}$.

The effect $f$ has on the coderivation $\hat{l}=\sum_n\hat{l}_n$ is this: $\hat{l}\mapsto f\circ\hat{l}\circ f^{-1}$.
To compute the difference, $\hat{l}-f\circ\hat{l}\circ f^{-1}$, remove the hats by defining 
$\beta_n=s\circ\hat{\beta}_n\circ((s^{-1})^{\otimes n})$, cf. sect.~\ref {def-via-coderiv}. As in the cited section, $\beta_n\in\Omega^n_A[1-n]$.
A somewhat painful but straightforward computation will then reveal that 
\[
\hat{l}-f\circ\hat{l}\circ f^{-1}= (d_{DR}+\text{Lie}_\partial)(\{\beta_n\})\hat{ },
\]
as desired. $\qed$

\section{beilinson-drinfeld}
\label{beilinson-drinfeld}

\subsection{ }
\label{from-n-to-any-set}

Let $\CV ect$ be the category {\em super}-vector spaces over $\BC$. For any collection $M_1,M_2,....,M_n\in\CV ect$ and any permutation $\sigma\in S_n$, there is a standard
isomorphism $\sigma: M_1\otimes M_2\otimes\cdots\otimes M_n\rightarrow M_{\sigma_1}\otimes M_{\sigma_2}\otimes\cdots\otimes M_{\sigma_n}$. For example,
$\sigma: M_1\otimes M_2\iso M_2\otimes M_1$
is defined by $\sigma(u\otimes v)= (-1)^{uv}v\otimes u$. Such isomorphisms satisfy an obvious associativity condition and make sense, for any finite set $I$, of
the tensor product $\otimes_I M_i$ of any $I$-family of vector spaces $\{M_i,\;i\in I\}$.

Suppose given a commutative, associative, purely even $\BC$-algebra $R$. If each $M_i$ above is an $R$-module, then we regard $M_1\otimes M_2\otimes\cdots\otimes M_n$
as an $R^n\stackrel{\text{def}}{=} R^{\otimes n}$-module. The tensor product $M_{\sigma_1}\otimes M_{\sigma_2}\otimes\cdots\otimes M_{\sigma_n}$ is also an $R^n$-module,
and we let $\sigma^*(M_{\sigma_1}\otimes M_{\sigma_2}\otimes\cdots\otimes M_{\sigma_n})$ to be an $R^n$-module that is a pull-back w.r.t. the algebra
isomorphism $\sigma: R^n\rightarrow R^n$. Then the map of vector spaces $\sigma$ becomes an $R^n$-module isomorphism
\[
\sigma: M_1\otimes M_2\otimes\cdots\otimes M_n\longrightarrow \sigma^*(M_{\sigma_1}\otimes M_{\sigma_2}\otimes\cdots\otimes M_{\sigma_n}).
\]
This allows one unambiguously to talk about an $R^I$-module $\otimes_I M_i$.

To push this a step further, assume now that $R$ is a  commutative Hopf  algebra with comultiplication $\Delta: R\rightarrow R^2$.   Iterating $\Delta$ one 
obtains an algebra morphism $R\rightarrow R^n$ for any $n$. There is a map
\[
Hom_{R^n}(\otimes_{i=1}^n M_{\sigma_i},N\otimes_R R^n)\longrightarrow Hom_{R^n}(\otimes_{i=1}^n M_i,N\otimes_R R^n) ,\; \phi\mapsto \phi^\sigma,
\]
$\phi^\sigma$ being defined as a unique morphism that makes the following diagram commutative
\begin{equation}
\label{comm-diagr-2}
\begin{gathered}
\xymatrix{ \otimes_{i=1}^nM_i\ar[r]^{\phi^\sigma}\ar[d]^{\sigma}&
N\otimes_R R^n\ar[d]^{1\otimes\sigma}\\
\otimes_{i=1}^nM_{\sigma_i}\ar[r]^{\phi}&N\otimes_R R^n}
\end{gathered}
 \end{equation}
 This makes sense of the space  $Hom_{R^I}(\otimes_{I} M_i,N\otimes_R R^I)$ for any set $I$.

\subsection{ }
\label{tansl-invar-version}

Now let $R=\BC[\partial]$ be a polinomial ring regarded as a Hopf algebra with comultiplication $\Delta: R\rightarrow R\otimes R$, $\partial\mapsto\partial\otimes 1+1\otimes\partial$,
and counit $\epsilon: R\rightarrow \BC$, $\partial\mapsto 0$, and
 let $\CM$ be the category of $R$-modules.

For a finite set $I$ and a collection of $R$-modules $M_i$, $i\in I$ and  $N$, define
\[
P^*_I(\{M_i\},N)= Hom_{R^I}(\otimes_{i\in I}M_i, N\otimes_R R^I).
\]

Elements of $P^*_I(\{M_i\},N)$ are called *-operations.

The composition is defined as follows:  for a surjection $\pi: J\twoheadrightarrow I$, and a collection of operations
$\psi_i\in P^*_{J_i}(\{L_j\}, M_i)$, $i\in I$, $J_i=\pi^{-1}(i)$, and $\phi\in P^*_I(\{M_i\},N)$, define $\phi(\psi_i)\in P^*_J(\{L_j\},N)$ to be the composite map
\footnote{To be precise, it is clear what this composition means if $J$ and $I$ are ordered in a way compatible with $\pi$; a change of orders leads to a change
of the composition map, which is easily shown to agree with the identifications made in sect.~\ref{from-n-to-any-set}. The interested reader is advised to read an
appendix to the book \cite{Lei} for more details.}:
\[
\otimes_JL_j\iso\otimes_I\otimes_{J_i}L_j\stackrel{\otimes\psi_i}{\longrightarrow}\otimes_I (M_i\otimes_R R^{J_i})=(\otimes_I M_i)\otimes_{R^I}R^J
\stackrel{\phi\otimes Id}{\longrightarrow} N\otimes_R R^I\otimes_{R^I}R^J=N\otimes_R R^J,
\]
where the identification $\otimes_I (M_i\otimes_R R^{J_i})=(\otimes_I M_i)\otimes_{R^I}R^J$ uses the natural algebra homomorphism
$R^I\rightarrow R^J$ that is the tensor product, over $I$, of the homomorphisms $R\rightarrow R^{J_i}$ obtained by iterating the coalgebra
map $R\rightarrow R^2$.

An associativity property holds: if, in addition, there is a surjection $K\twoheadrightarrow J$ and operations $\chi_j\in P^*_{H_j}(\{ A_k\},L_j)$, $j\in J$, then
$(\phi(\psi_i))(\chi)=\phi(\psi_i(\chi_j))$.  This defines on $\CM$  a
 {\em pseudo-tensor category} structure, \cite{BD}, 1.1.1\footnote{also
 called a symmetric multicategory in \cite{Lei}.}; $\CM$ equipped with this structure will be  denoted by $\CM^*$.  
 
We shall often encounter the situation when the $I$-family is constant, $M_i=M$, $J=I$ and $\pi$ is a bijection. In this case, the composition $\phi(id_M,id_M,...)$ also
belongs to $P^*_I(\{M\},N)$ and will be denoted $\phi^\pi$.  If $I=[n]=\{1,2,...,n\}$, then this defines a {\em left}  action of the symmetric group
\begin{equation}
\label{conve-fo-permu}
S_n\times P^*_{[n]}(\{M\},N)\longrightarrow P^*_{[n]}(\{M\},N);\; (\sigma,\phi)\mapsto \phi^\sigma, \phi^\sigma\stackrel{\text{def}}{=}\sigma^{-1}\circ\phi\circ\sigma,
\end{equation}
cf. (\ref{comm-diagr-2}); it would be correct if not pedantic to write $1\otimes \sigma^{-1}$ instead of $\sigma^{-1}$ in this formula.

For example, if $\phi(a,b)=\sum_{i,j}\langle a,b\rangle_{ij}\partial_1^i\partial_2^j$ for some  $\langle a,b\rangle_{ij}\in N$
and $\sigma=(1,2)$ is the transposition,
then
\[
\phi^\sigma(a,b)=\sum_{i,j}\langle b,a\rangle_{ij}\partial_2^i\partial_1^j=(1,2)\phi(b,a),
\]
provided both $a$ and $b$ are even.

More generally,  if 
\[
\phi(a_1,a_2,\ldots, a_n)=\sum_{i_1,\ldots,i_n}\langle a_1,a_2,\ldots, a_n\rangle_{i_1,\ldots,i_n}\partial_1^{i_1}\partial_2^{i_2}\cdots\partial_n^{i_n},
\]
then
\[
\phi^\sigma(a_1,a_2,\ldots, a_n)=\sum_{i_1,\ldots,i_n}\langle a_{\sigma_1},a_{\sigma_2},\ldots, a_{\sigma_n}\rangle_{i_1,\ldots,i_n}\partial_{\sigma_1}^{i_1}\partial_{\sigma_2}^{i_2}\cdots\partial_{\sigma_n}^{i_n},
\]
In both these formulas $\partial_i$ stands for $\underbrace{1\otimes\cdots\otimes\partial}_i\otimes\cdots\otimes 1\in R^n$.

\subsection{ }
\label{some-expl-form}
If a choice is made, then explicit formulas can be written down. If $I=\{1,2,3,...,n\}$, then $N\otimes_R R^I$ can be identified with $N[\partial_1,\partial_2,...,\partial_{n-1}]$. A binary operation $\mu\in P^*_{\{1,2\}}(\{M,M\},N)$  can then be written as follows
\[
\mu(a,b)=\sum_na_{(n)}b\otimes\frac{\partial_1^n}{n!}\text{ for some }a_{(n)}b\in N.
\]
One has for the transposition $\sigma=(1,2)$
\begin{eqnarray}
\mu^\sigma(a,b)&=&\sum_nb_{(n)}a\otimes\frac{\partial_2^n}{n!}=\sum_nb_{(n)}a\otimes\frac{(\partial_1+\partial_2-\partial_1)^n}{n!}\nonumber\\
&=&\sum_{n\geq j}(-1)^j(b_{(n)}a)\frac{\partial^{n-j}}{(n-j)!}\otimes\frac{\partial_1^j}{j!}\nonumber
\end{eqnarray}

Similarly, 
\[
\mu(a,\mu(b,c))=\sum_{n}\mu(a,b_{(n)}c)\otimes\frac{\partial_2^n}{n!}=\sum_{m,n}a_{(m)}(b_{(n)}c)\otimes\frac{\partial_1^m}{m!}\frac{\partial_2^n}{n!},
\]
but
\begin{eqnarray}
\mu(\mu(a,b),c)&=&\sum_{m}\mu(a_{(m)}b, c)\otimes\frac{\partial_1^m}{m!}=\sum_{m,n}(a_{(m)}b)_{(n)}c\otimes\frac{\partial_1^m}{m!}\frac{\Delta(\partial)^n}{n!}\nonumber\\
&=&\sum_{m,n}(a_{(m)}b)_{(n)}c\otimes\frac{\partial_1^m}{m!}\frac{(\partial_1+\partial_2)^n}{n!}\nonumber
\end{eqnarray}

\subsection{ }
\label{from pseudo-tensor-to tensor}
 A particular class of pseudo-tensor categories consists of tensor categories, where we use the words tensor
category as  synonymous with symmetric monoidal category. Namely, given such a category with $\otimes$ standing for the tensor product one defines
$P^*_I(\{M_i\},L)$ to be the Hom-space $Hom(\otimes_IM_i,L)$. Of course, not any pseudo-tensor category $\CC$ is a tensor category; to see the difference one
may want to consider the functor $\CC\rightarrow \CS ets$, $L\mapsto P^*_I(\{M_i\},L)$ and ask whether it is representable for each collection $\{M_i\}_I$.

The computations in \ref{some-expl-form} were intended to demonstrate how one works with a genuine pseudo-tensor category, but
it is sometimes useful to have a given pseudo-tensor category ``embedded'' as a subcategory of a tensor category. One universal such construction is suggested in \cite{BD}, 
Remark 1.1.6: Given a pseudo-tensor category $\CC$, define $\CC^\otimes$ as follows: an object is an  $I$-family  $\{M_i\}_I$, $M_i\in\CC$, a morphism 
$\{L_j\}_J\rightarrow \{M_i\}_I$ is a collection $\{\pi:J\twoheadrightarrow I, \phi_i\in P_{\pi^{-1}(i)}(\{L_j\},M_i)\}$. 

In the case of  our $\CM^*$, sect.~\ref{tansl-invar-version}, there is a more useful construction, \cite{BD}, 3.4.10.  In order to recall it, we need to introduce a bit of notation: given 2 sets $I$ and $J$, a surjection $\pi:J\twoheadrightarrow I$,
and an $R^I$-module $M$, let  $\Delta^{(\pi)}_*(M)$ be defined by
\[
\Delta^{(\pi)}_*M=M\otimes_{R^I}R^J,
\]
where the ring  morphism $R^I\rightarrow R^J$ is the one utilized in sect.~\ref{tansl-invar-version}.  The same object will sometimes be denoted by
$\Delta^{(J/I)}_*M$. If $I$ is a 1-element set, then we write simply  $\Delta^{(J)}_*M$.  For example, we have
\[
P^*_J(\{M_j\},L)=Hom_{R^I}(\otimes_I M_i,\Delta^{(J)}_*L).
\]
Note that $\Delta^{(\pi)}_*$ becomes a functor $Mod(R^I)\rightarrow Mod(R^J)$ if we define
\[
\Delta^{(\pi)}_*(f)\stackrel{\text{def}}{=}f\otimes Id:\; \Delta^{(\pi)}_*M\longrightarrow \Delta^{(\pi)}_*N\text{ for any  } f\in Hom_{R^I}(M,N).
\]
Similarly, given a diagram $K\stackrel{f}{\twoheadrightarrow} J \stackrel{g}{\twoheadrightarrow} I$, one has $\Delta^{(g\circ f)}_*=\Delta^{(f)}_*\circ\Delta^{(g)}_*$ or, if one wishes,
$\Delta^{(K/I)}_*=\Delta^{(K/J)}_*\circ\Delta^{(J/I)}_*$.

Now, the construction:

\subsubsection{ }
\label{catego-S}
{\em Define} an $R^\CS$-module $M$ ($\CS$ is to be thought of as ``Sets'') to be a rule that to each nonempty set $I$ assigns an $R^I$ module $M_{(I)}$ and to each surjection
$\pi:J\twoheadrightarrow I$ an $R^J$-module morphism $\theta^{(\pi)}_M: \Delta^{(\pi)}_*M_{(I)}\rightarrow M_{(J)}$ so that the following compatibility condition holds:
\[
\theta^{(f\circ g)}_M=\theta^{(g)}_M\circ\Delta^{(g)}_*(\theta^{(f)}_M).
\]
In what follows we will sometimes (for typographical reasons) suppress the index and write $\theta^{(\pi)}$ instead of $\theta_M^{(\pi)}$. Notations such as $\theta^{(J/I)}$ or
$\theta^{(J)}$ if $|I|=1$ will also be used.

Here is an example of an $R^\CS$-module: if $M$ is an $R$-module then define $\Delta^{(\CS)}_*M$ s.t.  $\Delta^{(\CS)}_*M_{(I)}=\Delta^{(I)}_*M$ with
$\theta^{(\pi)}_{\Delta^{(\CS)}_*M}=Id$. 

A morphism $f: M\rightarrow N$ is a collection of morphisms $f_{(I)}\in Hom_{R^I}(M_{(I)},N_{(I)})$ for each finite set $I$  that respects the structure, i.e.,  such that
\begin{equation}
\label{morph-in-tens-s-comp-cond}
\theta^{(\pi)}_N\circ\Delta^{(\pi)}_*(f_{(I)})=f_{(J)}\circ\theta^{(\pi)}_M.
\end{equation}
Denote by $Mod(R^\CS)$ the category of $R^\CS$-modules.  The construction of $\Delta^{(\CS)}_*M$ given above defines a functor:
\[
\Delta^{(\CS)}_*: Mod(R)\longrightarrow Mod(R^\CS),\; M\mapsto \Delta^{(\CS)}_*M,
\]
if we let $\Delta^{(\CS)}_*(f)_{(I)}=\Delta^{(I)}_*(f)$.  In fact, this functor is fully faithful, which follows from the following slightly more general observation (to be used more than once): 

for any $N\in R^\CS$ we have an isomorphism
\begin{equation}
\label{exte-morphi-2}
\Delta^{(\CS)}_*:\; Hom_{R}(M,N_{([1])})\iso Hom_{R^\CS}(\Delta^{(\CS)}_*M,N),\; f\mapsto\Delta^{(\CS)}_*(f)\text{ s.t. }\Delta^{(\CS)}(f)_{(J)}=\theta^{(J)}_N\circ\Delta^{(J)}_*(f).
\end{equation}
The fact that the the collection $\{\Delta^{(\CS)}(f)_{(J)}\}$ is a morphism is obvious. The map in the opposite direction 
$Hom_{R^\CS}(\Delta^{(\CS)}_*M,N)\rightarrow Hom_{R}(M,N_{[1]})$  is defined by $f\mapsto f_{([1])}$. The two maps are each other's
inverses because the definition of $\Delta^{(\CS)}_*(f)$ is forced on us by compatibility (\ref{morph-in-tens-s-comp-cond}).

We will sometimes take the liberty of informally referring to this phenomenon by saying that $\Delta^{(\CS)}_*M$ is freely generated by $M$.

\subsubsection{ }
\label{tnso-prduc}
Given an $I$-family of objects $\{M_i\}_I$ of   $Mod(R^\CS)$ define a {\em tensor product} $\otimes_I M_i$ by declaring:
\[
(\otimes_I M_i)_{(J)}=\bigoplus_{\pi:J\twoheadrightarrow I}\otimes_I(M_i)_{(\pi^{-1}(i))}.
\]
The structure morphisms $\theta^\alpha$, $\alpha:K\twoheadrightarrow J$, are naturally  defined:  notice that for each $\pi:J\twoheadrightarrow I$
\[
\Delta^{(\alpha)}_*(\otimes_I(M_i)_{\pi^{-1}(i)})=\otimes_I \Delta_*^{((\pi\circ\alpha)^{-1}(i)/\pi^{-1}(i))}(M_i)_{(\pi^{-1}(i))}.
\]
The tensor product of the structure morphisms $\otimes_I\theta_{M_i}^{\alpha_i}$, where $\alpha_i$ means the restriction of $\alpha$ to $(\pi\circ\alpha)^{-1}(i)$, gives then
the map 
\[
\Delta^{(\alpha)}_*(\otimes_I(M_i)_{(\pi^{-1}(i))})\longrightarrow \otimes_I (M_i)_{(\pi\circ\alpha)^{-1}(i))}.
\]
Summation over $\pi$ gives a map
\[
\Delta^{(\alpha)}_*(\otimes_I M_i)_{(J)}\longrightarrow \bigoplus_{\pi:J\twoheadrightarrow I}\otimes_I(M_i)_{((\pi\circ\alpha)^{-1}(i))}.
\]
The target of this map  is clearly a subset of 
of $(\otimes_I M_i)_{(K)}$ -- the one that involves only those surjections $K\twoheadrightarrow I$ that factor through $\alpha$; the embedding of the subset into the set is our 
$\theta^{(\alpha)}$.

In particular, if $M_i=\Delta^{(\CS)}_*N_i$, then the tensor product  is given by
\begin{equation}
\label{defn-tens-in-s-of-from-m-star}
(\otimes_I \Delta^{(\CS)}_*N_i)_{(J)}=\bigoplus _{J\twoheadrightarrow I}\Delta^{(J/I)}_* (\otimes_I N_i).
\end{equation}
For example, if $J=I$, then this becomes a sum over   bijections of $I$ on itself
\begin{equation}
\label{defn-tens-in-s-of-from-m-star-part-case}
(\otimes_I \Delta^{(\CS)}_*N_i)_{(I)}=\bigoplus _{\sigma\in\text{Bij}(I)}\Delta^{(\sigma)}_*( \otimes_I N_i),
\end{equation}
and not just $\otimes_I N_i$, as one could naively expect (at which point the reader is invited to figure out the meaning of $\Delta^{(\sigma)}_*( \otimes_I N_i)$.)

The meaning of (\ref{defn-tens-in-s-of-from-m-star}) is that $\otimes_I \Delta^{(\CS)}_*N_i$ is freely generated by $\otimes_I N_i$.
Namely, there is an isomorphism, cf. (\ref{exte-morphi-2}),
\begin{equation}
\label{exte-morphi-3}
\Delta^{(\CS)}_*:\; Hom_{R^I}(\otimes_I N_i,L_{(I)})\iso Hom_{R^\CS}(\otimes_I \Delta^{(\CS)}_*N_i, L),
\end{equation}
where $\Delta^{(\CS)}_*(f)$ is defined by the familiar requirement:  $\Delta^{(\CS)}_*(f)_{(J)}$ restricted to $\Delta^{(J/I)}_*\otimes_I N_i$ 
equals $\theta^{(J)}_N\circ\Delta^{(J/I)}_*(f)$.

Indeed, the map in the opposite direction $Hom_{R^\CS}(\otimes_I \Delta^{(\CS)}_*N_i, L)\rightarrow Hom_{R^I}(\otimes_I N_i,L_{(I)})$  is defined by 
$f\mapsto f_{(I)}|_{\otimes_I N_i}$, and the rest of the proof is as that of (\ref{exte-morphi-2}).
\subsubsection{ }
\label{pse-te-fctr}
The reason this discussion has been undertaken is the following result: the fully faithful functor $\Delta^{(\CS)}_*: Mod(R)\longrightarrow Mod(R^\CS),\; M\mapsto \Delta^{(\CS)}_*M$ introduced  above is, in fact, a fully faithful pseudo-tensor  functor $\Delta^{(\CS)}_*: \CM^*\longrightarrow Mod(R^\CS),\; M\mapsto \Delta^{(\CS)}_*M$, i.e., there is a natural vector space isomorphism
\begin{equation}
P^*_I(\{M_i\},L)\iso \text{Hom}_{\text{Mod}(R^\CS)}(\otimes_I \Delta^{(\CS)}_*M_i, \Delta^{(\CS)}_*L).
\label{full-fai-pr-tens-fnctr}
\end{equation}
This is a particular case of (\ref{exte-morphi-3}).

\subsection{ }
\label{augm-fnctr}
Along with $\CM^*$ consider $\CV ect$, the tensor category of vector spaces, hence a pseudo-tensor category where $P_I(\{V_i\}, V)=Hom_\BC(\otimes_I V_i,V)$.
The assignment $\CM^*\mapsto h(M)\stackrel{\text{def}}{=}M/\partial M$ defines a pseudo-tensor functor, called an {\em augmentation functor} in \cite{BD}, 1.2.4, 1.2.9-11,
\[
h:\CM^*\longrightarrow \CV ect,
\]
as $h$   defines, in an obvious manner, a map
\[
h_I:\; P^*_I(\{M_i\},N)\longrightarrow Hom_\BC(\otimes_I h(M_i),h(N)),
\]
which is functorial in $\{M_i\}$ and $N$.

\subsection{ }
\label{dens-star-lie-ass-etc}
A pseudo-tensor category structure, i.e., a family of well-behaved spaces of ``operations'' $P^*_I(\{L_i\},N)$, is what is needed to define various algebraic structures.
For example,  a Lie* or associative* algebra is a pseudo-tensor functor
\[
Lie\longrightarrow \CM^*\text { or } Ass\longrightarrow \CM^*,
\]
where $Lie$ or $Ass$ (resp.) is the corresponding operad (an operad being a pseudo-tensor category with a single object.)  Explicitly, this means a choice of  an $R$-module $V$ and  an operation $\mu(.,.)\in P^*_{\{1,2\}}(\{V,V\},V)$
that satisfies appropriate identities written by means of the above defined composition. For example, $V$ is an associative*  algebra  if
$\mu(\mu(.,.),.)=\mu(.,\mu(.,.))$ as elements of $P^*_{\{1,2,3\}}(\{V,V,V\},V)$. Likewise, $V$ is Lie* if 
\[
\mu^{(1,2)}(.,.)=-\mu(.,.)) \text{ and } \mu(\mu(.,.),.) + \mu(\mu(.,.),.)^{(1,2,3)}+\mu(\mu(.,.),.)^{(1,2,3)^2}=0,
\]
see (\ref{conve-fo-permu}) for some of the notation used.

It is easy to verify, using \ref{some-expl-form},  that  a Lie* algebra is an $R$-module $V$ with a family of multiplications $_{(n)}$ s.t. $a_{(n)}b=0$ if $n\gg 0$ and (assuming for simplicity  that $V$ is purely even)
\begin{eqnarray}
a_{(n)}b=(-1)^{n+1}\sum_{j\geq 0}(b_{(n+j)}a)\frac{\partial^j}{j!}& &\text{: anti-symmetry}\label{antisymmetry}\\
 a_{(n)}b_{(m)}c-b_{(m)}a_{(n)}c=\sum_{j\geq 0}{n\choose j}(a_{(j)}b)_{(n+m-j)}c& &\text{: Jacobi}\label{jacobi}
\end{eqnarray}

The last equality is known as the {\em Borcherds commutator formula}

It is convenient to denote by $a(\partial)$ the formal sum $\sum_n a_{(n)}\partial^n/n!$.
We have

(i)  the just written Jacobi identity is equivalent to
\[
a(\partial_1)b(\partial_2)c-b(\partial_2)a(\partial_1)c=(a(\partial_1)b)(\partial_1+\partial_2)c.
\]
(ii)  the  associativity condition $\mu(\mu(.,.),.)=\mu(.,\mu(.,.))$ is equivalent to
\[
a(\partial_1)b(\partial_2)c=(a(\partial_1)b)(\partial_1+\partial_2)c.
\]
This point of view has been introduced and developed by V.Kac and his collaborators, see \cite{K} and references therein, especially \cite{BKV}, sect.~12.

\subsection{ }
\label{lie-star-chevalley}
Let $L$ be a Lie* algebra with bracket $[.,.]\in P^*_{\{1,2\}}(\{L,L\},L)$.  An $R$-module $M$ is called an $L$-module if
there is an operation $\mu\in P^*_{\{1,2\}}(\{L,M\},M)$ s.t. 
\[
\mu(.,\mu(.,.))-\mu(.,\mu(.,.))^{(1,2)}=\mu([.,.],.).
\]
The untiring reader will have no trouble verifying that in terms of $_{(n)}$-products this is nothing but an obvious version
of (\ref{jacobi}).

The Chevalley complex is defined as follows. Denote by $C^n(L,M)$ the subspace of $P^*_{[n]}(\{L\},M)$, $[n]=\{1,2,...,n\}$, of anti-invariants of the symmetric group action.  
Set, mimicking the usual definition,
\begin{eqnarray}
& &d: C^n(L,M)\longrightarrow C^{n+1}(L,M)\text{ s.t.}\nonumber\\
 & &d\phi(l_1,l_2,\ldots,l_{n+1})=
 \sum_{1\leq i\leq n+1}(-1)^{i+1}\sigma^{-1}_i\circ\mu(l_i,\phi(l_1,...,\widehat{l_i},...,l_{n+1}))+\nonumber\\
 & &\sum_{1\leq i< j\leq n+1}(-1)^{i+j}\sigma^{-1}_{ij}\circ\phi([l_i,l_j],l_1,...,\widehat{l_i},...
 ,\widehat{l_j},....,l_{n+1})^{}.\nonumber
\end{eqnarray}
where $\sigma_i$  ($ \sigma_{ij}$ ) stands for  the permutation applied to the variables of the first ( last)  term; cf. (\ref{conve-fo-permu}). Essentially the familiar (from ordinary Lie theory) proof shows that $d^2=0$.

Various computations involving this complex, called there {\em reduced},  can be found in \cite{BKV}.

\subsection{ }
\label{from lie*tolie}
If $L$ is a Lie* algebra and $M$ an $L$-module, then $h(L)$ is an ordinary Lie algebra and $h(M)$, as well as $M$ itself is an $h(L)$-module. This is true on general grounds, see  sect.~\ref{augm-fnctr}, but also easily follows from the explicit formulas of sect.~\ref{dens-star-lie-ass-etc}.

\subsection{ }
\label{compound}
In order to define a Poisson algebra object in $\CM^*$ one needs, in addition to Lie*, another structure, associative commutative multiplication, and
another constraint, the Leibniz rule. This is taken care of by another pseudo-tensor structure on $\CM$, in fact, a genuine tensor category structure engendered
by the fact that $R$ is a Hopf algebra.  Given $A,B\in\CM$, let $A\shriekotimes B$  be $A\otimes B$ acted upon by $R$ via $\Delta: R\rightarrow R\otimes R$.
The category $\CM$ with this tensor structure will be denoted by $\CM^!$.

The 2 pseudo-tensor structures are related in that operations can sometimes be multiplied. Let us describe this product in the simplest possible case. 
Assume given $P^*_I(\{M_i\},N_1)$,  $P^*_J(\{L_j\},N_1)$, and fix $i_0\in I$, $j_0\in J$. Denote by  $I\vee J$ (or rather $I\vee_{i_0j_0}J$)
the  union $I\sqcup J$ modulo the relation $i_0=j_0$. There is a natural map
\begin{equation}
\label{defn-shriek-tens}
P^*_I(\{M_i\},A)\otimes P^*_J(\{L_j\},B)\longrightarrow P^*_{I\vee J}(\{M_i,L_j,M_{i_0}\shriekotimes L_{j_0}\}_{i\not=i_0,j\not=j_0},A\shriekotimes B).
\end{equation}

In order to define this map, it is best to build up on sect.~\ref{from pseudo-tensor-to tensor} and introduce some pull-back functors
operating among the categories of $R^I$-modules. Namely, given a surjection $\pi: J\twoheadrightarrow I$  an $R^I$-module $M$, and
an $R^J$-module $N$
in addition to
to $\Delta^{(J/I)}_*M=M\otimes_{R^I}R^J\in Mod(R^J)$, see {\em loc. cit.}, which has the meaning of a push-forward, define 
$\Delta^{(J/I)*}N\in Mod(R^I)$ to be the pull-back w.r.t. $R^I\rightarrow R^J$. One has  $A\shriekotimes B=\Delta^{([2])*}A\otimes B$.
More generally, there is an obvious projection $I\sqcup J\twoheadrightarrow I\vee J$ and an isomorphism of  $R^{I\vee J}$-modules
\[
\otimes_{i\not=i_0} M_i\otimes\otimes_{j\not=j_0} L_j\otimes(M_{i_0}\shriekotimes N_{j_0})\iso
\Delta^{(I\sqcup J/I\vee J)*}(\otimes_IM_i\otimes\otimes_J L_j).
\]
This implies that given
\[
\phi\in P^*_I(\{M_i\},A)=Hom_{R^I}(\otimes_I M_i,\Delta^{(I)}_*A),\;
\psi\in P^*_J(\{L_j\},B)=Hom_{R^J}(\otimes_J L_j,\Delta^{(J)}_*B),
\]
we obtain
\[
\phi\otimes\psi\in Hom_{R^{I\sqcup J}}(\otimes_I M\otimes\otimes_J L_j,\Delta^{(I)}_*A\otimes \Delta^{(J)}_*B)=
Hom_{R^{I\sqcup J}}(\otimes_I M\otimes\otimes_J L_j,\Delta^{(I\sqcup J/[2])}_*A\otimes B),
\]
and pulling back
\[
\Delta^{(I\sqcup J/I\vee J)*}(\phi\otimes\psi)\in
Hom_{R^{I\vee J}}(\Delta^{(I\sqcup J/I\vee J)*}(\otimes_I M\otimes\otimes_J L_j),
\Delta^{(I\sqcup J/I\vee J)*}\Delta^{(I\sqcup J/[2])}_*A\otimes B).
\]
There is a base change isomorphism of functors
\[
\Phi_\bullet: \Delta^{(I\sqcup J/I\vee J)*}\Delta^{(I\sqcup J/[2])}_*(\bullet)\iso\Delta^{(I\vee J)}_*\Delta^{([2])*}(\bullet).
\]
Indeed, for the former we have
\[
\Delta^{(I\sqcup J/I\vee J)*}\Delta^{(I\sqcup J/[2])}_*(A\otimes B)=\Delta^{(I\sqcup J/I\vee J)*}(A\otimes B)\otimes_{R^2} R^{I\sqcup J},
\]
for the latter
\[
\Delta^{(I\vee J)}_*\Delta^{([2])*}(A\otimes B)=(A\shriekotimes B)\otimes_R R^{I\vee J}=(A\otimes B)\otimes_{R^2}({R^2}\otimes_R R^{I\vee J}).
\]
Consider the ring morphism
\begin{equation}
\label{defnofshriek}
R^2\otimes_R R^{I\vee J}\longrightarrow R^I\otimes R^J
\end{equation}
defined on the generators to be the following two:
\[
R^2\longrightarrow  R^I\otimes R^J\text{ and } R^{I\vee J}\longrightarrow R^I\otimes R^J.
\]
The former is the tensor product of the iterated coproduct maps $R\rightarrow R^I$ and $R\rightarrow R^J$. The latter is defined to be
$\partial_\alpha\mapsto\partial_\alpha$ if $\alpha$ is different from the equivalence class $\{i_0,j_0\}$ and 
$\partial_\alpha\mapsto\partial_{i_0}+\partial_{j_0}$ if $\alpha$ is the equivalence class $\{i_0,j_0\}$. The map (\ref{defnofshriek}) is an isomorphism
as it is simply a coordinate  change in a polynomial ring. $\Phi_{A\otimes B}$ is induced by the inverse of (\ref{defnofshriek}).

Product (\ref{defn-shriek-tens}) is defined as follows:
\[
\phi\otimes\psi\mapsto \Phi_{A\otimes B}\circ \Delta^{(I\sqcup J/I\vee J)*}(\phi\otimes\psi).
\]

 Denote by $\phi\shriekotimes_{i_0,j_0}\psi$ the tensor product of 2 operations thus defined.

\subsection{ }
\label{exs-of-prod-oper}
\begin{sloppypar}
If index sets are ordered and operations are written in terms of $_{(n)}$-products, sect.~\ref{some-expl-form}, then the inherent symmetry of the definition is destroyed.
For example, given $\phi\in P^*_{\{1,2\}}(\{M_1,M_2\}, N)$, with $\phi(a,b)= \sum_n (a_{(n)}b)\otimes \frac{\partial_1^n}{n!}$ and $id\in P^*_{\{1\}}(L,L)$ one easily computes 
$\phi\shriekotimes_{21} id\in P^*_{\{1,\{2,1\}\}}(\{M_1,M_2\shriekotimes L\},N\shriekotimes L)$ to be
\end{sloppypar}
\begin{equation}
\label{shriek-tens-one}
\phi\shriekotimes_{21} id(a,b\shriekotimes c)=\sum_n (a_{(n)}b)\shriekotimes c\otimes \frac{\partial_1^n}{n!}.
\end{equation}
On the other hand,  $id\shriekotimes_{11} \phi\in P^*_{\{\{1,1\},2\}}(\{L\shriekotimes M_1,M_2\},L\shriekotimes N)$ is as follows
\begin{equation}
\label{shriek-tens-two}
id\shriekotimes_{11} \phi(c\shriekotimes a,b)=\sum_{n\geq j} (-1)^{n-j}(c\frac{\partial^{n-j}}{(n-j)!})\shriekotimes a_{(n)}b)\otimes \frac{\partial_{11}^j}{j!}.
\end{equation}
Indeed,  the construction of map (\ref{defn-shriek-tens}) gives the composition
\begin{multline}
c\shriekotimes a\otimes b\mapsto \sum_nc\otimes a_{(n)}b\otimes\frac{\partial_2^n}{n!}\mapsto\sum_nc\otimes a_{(n)}b\otimes\frac{(\partial_{11}-\partial_1)^n}{n!}\\
=\sum_{n\geq j}(-1)^{n-j}c\otimes a_{(n)}b\otimes\frac{\partial_{11}^j\partial_1^{n-j}}{j!(n-j)!}
=\sum_{n\geq j}(-1)^{n-j}(c\frac{\partial^{n-j}}{(n-j)!})\shriekotimes a_{(n)}b\otimes\frac{\partial_{11}^j}{j!}\nonumber,
\end{multline}
as desired.  In this computation,  the last equality follows from the fact that
$id\in P^*_{\{1\}}(L,L)= Hom_R(L,L\otimes_RR)$, and so  $c\otimes \partial_1^m= c\partial^m\otimes 1$.

\subsection{ }
\label{comm-modul-act-prod-comm}
A commutative$^{!}$ algebra is defined to be a commutative (associative unital ) algebra in $\CM^!$.  In the present context, this is the same thing as the
conventional commutative (associative unital) algebra with derivation. Modules over a commutative$^!$ algebra are defined (and described) similarly.

If  $A$ is a commutative and associative algebra in $\CM^!$, $\{M_i\}_{[n]}$ and $N$ are $A$-modules ,  with actions $\mu_i: A\otimes M_i\rightarrow M_i$ and
$\mu:A\otimes N\rightarrow N$, then an operation $\nu\in P^*_{[n]}(\{M_i\},N)$
is called $A$-multilinear if  
\begin{equation}
\label{defi-multil-in shrie-ca}
\mu(Id_A\otimes^!_{\{1,i\}}\nu)=\nu(Id,Id,...,\mu_i, Id,...,Id)
\end{equation}
 for all $i\in[n]$.

If $(L,[,.])$ is a Lie* algebra and $(M_1,\mu_1)$, $(M_2,\mu_2)$ are $L$-modules, $\mu_j\in P^*_{\{1,2\}}(\{L,M_j\},M_j)$ being the action,
$j=1,2$, then $M_1\shriekotimes M_2$ carries an $L$-module structure defined via the Leibniz
rule. Namely, one defines $\mu=\mu_1\shriekotimes id_{M_2}+\mu_2\shriekotimes id_{M_1}\in P^*_{\{1,2\}}(\{L,M_1\shriekotimes M_2\},M_1\shriekotimes M_2)$
and verifies, just as in the ordinary Lie algebra case, that this is a Lie* action.

If $A$ is a commutative$^!$ algebra, then we say that $L$ acts on $A$ (or $L$ acts on it by derivations) if  $A$ is an $L$-module s.t. the multiplication morphism
\[
A\shriekotimes A\longrightarrow A
\]
is a morphism of $L$-modules.

In a similar vein, $\CL$ is a Lie* $A$-algebroid if it is a Lie* algebra, an $A$-module, and it acts on $A$ (by derivations) s.t.

(1) the action $\mu\in P^*_{\{1,2\}}(\{\CL,A\},A)$ is $A$-linear w.r.t. the $\CL$-argument;

(2)  the $A$-module morphism
\[
A\shriekotimes\CL\longrightarrow\CL
\]
is an $\CL$-module morphism, cf. (\ref{lie-alg-comp-str}).

A coisson algebra $\CP$ is a Lie* algebra and a commutative$^!$ algebra s.t. the commutative$^!$-product map
\[
\CP\shriekotimes\CP\longrightarrow\CP
\]
is a Lie* algebra module morphism.

\subsection{ }
\label{ poiss-jet-coiss}
Let $A$ be a conventional commutative associative unital algebra. Denote by $J_\infty A$ the universal commutative associative algebra with derivation
generated by $A$. More formally, $J_\infty$ is the left adjoint of the forgetful functor from the category of commutative algebras with
derivation to the category of commutative algebras.
\begin{lem}.
\label{ poiss-jet-coiss-lem}
If $A$ is a Poisson algebra, then $J_\infty A$ is canonically a coisson algebra.
\end{lem}

{\em Proof.} If $\{.,.\}$ is the Poisson bracket on $A$, then define $\{a\partial^l,b\partial^k\}=\{a,b\}\otimes \partial_1^l\partial_2^k$, and extend to all of $J_\infty A$
using the Leibniz property; this makes perfect sense thanks to the universal property of $J_\infty A$. The relation 
$\{a,(bc)\partial\}=\{a,(b\partial)c+b(c\partial)\}$
is almost tautological. Indeed, by construction, the L.H.S. equals  $\{a,bc\}\otimes\partial_2$, the R.H.S. is  the sum of 4 terms
\[
\{a,b\partial\}c+\{a,c\}(b\partial)+\{a,b\}(c\partial)+\{a,c\partial\}b;
\]
the 1st plus the 3rd equals $\{a,b\}c\otimes\partial_2$, the 2nd plus the 4th equals $\{a,c\}b\otimes\partial_2$, which adds to  $\{a,bc\}\otimes\partial_2$, as desired.
$\qed$

In hindsight, this simple assertion appears to be this theory's {\em raison d'\^etre}.

To see an example, let $A$ be a commutative algebra and consider the symmetric algebra $S^\bullet _AT_A$, which is canonically Poisson, sect.~\ref{defn-of -tdo}.
It is graded, by assigning degree 1 to $T_A$, and so is the coisson algebra $J_\infty S^\bullet_AT_A$. Consider its degree 1 component, $J_\infty T_A$, which, by the way,
can be equivalently described as the universal $J_\infty A$-module with derivation generated by $T_A$. The Lie* bracket on $J_\infty S^\bullet_AT_A$ restricts to $J_\infty T_A$ and makes
it a Lie* algebra. Furthermore, $J_\infty S^\bullet_AT_A$ is a $J_\infty T_A$-module and $J_\infty A\subset J_\infty S^\bullet_AT_A$ is a submodule. Hence
$J_\infty T_A$ acts on $J_\infty A$ be derivations. One easily verifies that, in fact, $J_\infty T_A$ is a Lie* $J_\infty A$-algebroid, sect.~\ref{comm-modul-act-prod-comm}. Furthermore, it is not hard to prove that
if a Lie* algebra $L$ acts on $J_\infty A$ by derivations, then this action factors through a Lie* algebra morphism $L\rightarrow J_\infty T_A$.

A much more general discussion of tangent algebroids can be found in \cite{BD}, 1.4.16.

\subsection{ }
\label{lie-star-infty}
The context of Lie* brackets makes it straightforward to suggest a definition of a $\text{Lie}^*_\infty$ algebra; in what follows we will freely use the
notation of  sects.~\ref{lie-infty}, \ref{def-via-coderiv}.  

We shall say that an $R$-module $V$ is graded if $V=\oplus_{i\in\BZ} V_i$ and $R(V_i)\subset V_i$. Similarly, if $\{V_i\}$ and $W$ are graded
$R$-modules, we shall say that an operation $\mu\in P^*_I(\{V_i\},W)$ has degree $k$ if
\[
\mu(v_1,v_2\ldots)\in V_N\otimes_RR^I,\text{ where }\sum_i\text{deg}v_i+k=N.
\]
Similarly, if $V$ is a graded $R$-module, we shall say that an operation $\mu\in P^*_{[n]}(\{V\},W)$ is {\em antisymmetric} if 
\[
\mu^{\sigma}(v_{1},v_{2}\ldots v_{n})=\text{sgn}(\sigma)\epsilon(\sigma,\vec{x})\mu(v_1,v_2\ldots v_n),
\]
where the action of the symmetric group on operations, $(\sigma,\nu(...))\mapsto \nu^ \sigma(...)$ is the one defined in sect.~\ref{tansl-invar-version}.

Denote by  $\Lambda P^*_{[n]}(\{V\},W)$ the subspace of all antisymmetric operations.

{\em Definition.} A $\text{Lie}^*_\infty$ algebra is a graded $R$-module $L$ and a collection of antisymmetric $[n]$-operations
$l_n\in \Lambda P^*_{[n]}(\{L\},L)$, $\text{deg}l_n=2-n$, that for each $k=1,2\ldots$ satisfy the following identity
\begin{equation}
\label{jacobi-for-infty}
\sum_{i+j=k+1}\sum_{\sigma}\text{sgn}(\sigma)\epsilon(\sigma,\vec{x})(-1)^{i(j-1)}\sigma^{-1}\circ l_j(l_i(x_{\sigma_1}\ldots x_{\sigma_i}),x_{\sigma_{i+1}}\ldots x_{\sigma_n}),
\end{equation}
where $\sigma$ runs through the set of all $(i,n-i)$ unshuffles, i.e.,  $\sigma\in S_n$ s.t. 
$\sigma_1<\sigma_2<\cdots<\sigma_i$ and $\sigma_{i+1}<\sigma_{i+2}<\cdots<\sigma_n$; the meaning of  $\sigma\circ$ is as in (\ref{conve-fo-permu}).

\bigskip

By definition, $l_1$ is simply a degree 1 linear map $L\longrightarrow L$, and (\ref{jacobi-for-infty}) with $n=1$ says that $l_1^2=0$; in other words, $(L,l_1)$ is
a complex.

Let us denote $l_2(.,.)$ by $[.,.]$. One has $[x_1,x_2]=-(-1)^{x_1x_2}[x_2,x_1]^{(1,2)}$, and (\ref{jacobi-for-infty}) with $n=2$ reads, after an obvious re-arrangement,
\[
l_1[x_1,x_2]=[l_1(x_1),x_2]+(-1)^{x_1}[x_1,l(x_2)].
\]
We conclude that $[.,.]$ is an antisymmetric super-star-bracket of degree 0, and $l_1$ is its derivation. More explicitly, if we write $[x,y]=\sum_ix_{(i)}y\otimes \partial_1^i/i!$, then
\[
l_1(x_{(i)}y)=l_1(x)_{(i)}y+(-1)^xx_{(i)}l_1(y),
\]
hence $l_1$ is a derivation of all products $_{(i)}$.

The $n=3$ case of (\ref{jacobi-for-infty}) involves terms such as $[[.,.],.]$, $l_3\circ l_1$, and $l_1\circ l_3$. The first one will give the ``jacobiator,'' the last two will
show that the super-Jacobi identity holds up to homotopy, $l_3$:
\begin{multline}
[[x_1,x_2],x_3]+(-1)^{x_3(x_1+x_2)}[[x_3,x_1],x_2]^\sigma+(-1)^{x_1(x_2+x_3)}[[x_2,x_3],x_1]^{\sigma^2}=\\
-\left( l_1l_3(x_1,x_2,x_3)+l_3(l_1(x_1),x_2,x_3)+(-1)^{x_1}l_3(x_1,l_1(x_2),x_3) +(-1)^{x_1+x_2}l_3(x_1,x_2,l_1(x_3))\right)\nonumber
\end{multline}
Writing $l_3(x_1,x_2,x_3)=\sum_{m,n}(x_1,x_2,x_3)_{mn}\otimes\partial_1^m\partial_2^n/m!n!$ and equating the terms in front of $\partial_1^m\partial_2^n$ in the last
equality, we
obtain, cf. (\ref{jacobi}),
\begin{multline}
a_{(n)}b_{(m)}c-(-1)^{ab}b_{(m)}a_{(n)}c-\sum_{j\geq 0}{n\choose j}(a_{(j)}b)_{(n+m-j)}c=\\
l_1((a,b,c)_{mn})+(l_1(a),b,c)_{mn}+(-1)^a(a,l_1(b),c)_{mn}+(-1)^{a+b}(a,b,l_1(c))_{mn}\nonumber,
\end{multline}
where we took the liberty of using $a,b,c$ in place of $x_1,x_2,x_3$ (resp.) so as to avoid being flooded by  indices.

It is clear, of course, how the concept of a  differential Lie* superalgebras is defined and how that   of a $\text{Lie}^*_\infty$ algebra generalizes it.

\subsection{ }
\label{liestarinfty-via-coderiv}
In order to push the analogy with  ordinary $\text{Lie}_\infty$ algebras a little  further,  we would like to find an appropriate generalization of the material recalled in
 sect.~\ref{def-via-coderiv}, i.e., we seek a ``coalgebra with square 0 coderivation.'' The problem here is that if $L$ is a $\text{Lie}^*_\infty$ algebra, then
 we need an object of the type $S^\bullet L[1]= \oplus S^n L[1]$, but none of the symmetric powers $S^n L[1]$ is an object of $\CM^*$ unless $n=1$. More technically,
 the difficulty is that $\CM^*$ is not a tensor category, and this is where the construction of sect.~\ref{from pseudo-tensor-to tensor} is useful.
 
 \subsubsection{ }
 \label{star-comulti-tenso-case}
 As a warm-up, let us do the tensor algebra case, cf. Remark~\ref{remark-on-tens-vs-assinfty-alg}. Given an $R^\CS$-module $M$, sect.~\ref{catego-S},
  let $TM$ be the {\em reduced} free associative
 algebra generated by $M$ (in $Mod(R^\CS)$.)  This is nothing but  $\oplus_{n\geq 1}\otimes_{[n]}M$, where $[n]=\{1,2,,...,n\}$, see
 sect.~\ref{tnso-prduc} for the definition of $\otimes$; the word reduced  means that the algebra does not have a unit,
 a slight complication stemming from the fact that $\otimes_\emptyset M$ does not make much sense in our situation. 
 
 In the case where $M=\Delta^{(\CS)}_*L$, in addition to being an algebra, $T\Delta^{(\CS)}_*L$ carries a coalgebra structure. To define it, notice
 that, for each $n$, the ordinary comultiplication (\ref{defn-comul-us-te-alge}) defines a morphism
 $\Delta\in Hom_{R^{[n]}}(\otimes_{[n]}L, \oplus_{i\in [n]}(\otimes_{[i]}L)\otimes(\otimes_{[n-i]}L))$. 
 If we regard $\otimes_{[i]}L$ as $(\otimes_{[i]}\Delta^{(\CS)}_*L)_{([i])}$ and  
 $\otimes_{[n-i]}L$ as $(\otimes_{[n-i]}\Delta^{(\CS)}_*L)_{([n-i])}$, then this gives us, due to (\ref{exte-morphi-3}), a morphism
 \[
 \Delta^{(\CS)}_*(\Delta):\;\otimes_{[n]}\Delta^{(\CS)}_*L\longrightarrow  \oplus_{i\in [n]}(\otimes_{[i]}\Delta^{(\CS)}_*L)\otimes(\otimes_{[n-i]}\Delta^{(\CS)}_*L)).
 \]
 Hence a morphism, to be denoted in the same way,
 \[
 \Delta^{(\CS)}_*(\Delta):\; T\Delta^{(\CS)}_*L\longrightarrow T\Delta^{(\CS)}_*L\otimes T\Delta^{(\CS)}_*L.
 \]
 Explicitly,  
 \[
 \Delta_*^{(J/[n])}(\otimes_{[n]}L)\ni x_1x_2\cdots x_n\otimes P\mapsto \sum_{i=1}^{n-1}x_1x_2\cdots x_i\otimes x_{i+1}x_{i+2}\cdots x_n\otimes P,\; P\in R^J.
 \]
 This formula is not obviously different from (\ref{defn-comul-us-te-alge}), except that the terms of the type $1\otimes x$ or $x\otimes 1$ are missing, but notice a subtlety:  the notation $x_1x_2\cdots x_n\otimes P$ is ambiguous
 and its meaning depends on a choice of a surjection $J\twoheadrightarrow I$.
 
 The coassociativity is an immediate consequence of that of $\Delta$, and so $(T\Delta^{(\CS)}_*L,\Delta^{(\CS)}_*(\Delta))$ is a coalgebra in $R^\CS$.
 
 The concept of {\em coderivation} is defined in an obvious way to be a morphism $\delta:T\Delta^{(\CS)}_*L\rightarrow T\Delta^{(\CS)}_*L$  s.t.
 $(\delta\otimes 1+1\otimes\delta)\circ\Delta^{(\CS)}(\Delta)= \Delta^{(\CS)}(\Delta)\circ\delta$.   We shall now attach a coderivation to any morphism 
 $f\in Hom_{R^{[m]}}(\otimes_{[m]}L,\Delta^{([m])}_*L)$.

 Notice that given  $f\in Hom_{R^{[m]}}(\otimes_{[m]}L,\Delta^{([m])}_*L)$ and  numbers $i,n$ s.t. $1\leq i\leq n-m$, there arises a morphism $f_i$ defined to be the composition
 \begin{eqnarray}
 & &\otimes_{[n]}L=\otimes_{[1,i]}L\otimes\otimes_{[i+1,i+m]}L\otimes\otimes_{[i+m+1,n]}L\stackrel{1\otimes f\otimes 1}{\longrightarrow}
 \otimes_{[1,i]}L\otimes\Delta^{([i+1,i+m])}_*L\otimes\otimes_{[i+m+1,n]}L=\nonumber\\
 & &\Delta^{([n]/[n-m+1])}_*\otimes_{[n-m+1]}L.\nonumber
 \end{eqnarray}
 Then we obtain $\Delta^{(\CS)}_*(f)_i\in End_{R^{(\CS)}}(T\Delta^{(\CS)}_*L)$. Explicitly, this morphism operates according to a familiar-looking formula:
 \[
 x_1x_2\cdots x_n\otimes P\mapsto x_1\cdots x_{i}f(x_{i+1}\cdots x_{i+m})x_{i+m+1}\cdots x_n\otimes P,\; P\in R^n,
 \]
 but, of course, now $f(x_{i+1}\cdots x_{i+m})$ is not an element of $L$. Set
 \begin{equation}
 \label{const-code-tenso-chira}
 \delta(f)=\sum_{i=1}^\infty\Delta^{(\CS)}_*(f)_i.
 \end{equation}
 By construction, this is a well-defined element of $End_{R^{(\CS)}}(T\Delta^{(\CS)}_*L)$.  Furthermore, $\delta(f)(\otimes_{[n]}\Delta^{(\CS)}_*L)\subset
 \otimes_{[n-m+1]}\Delta^{(\CS)}_*L$, and so we can say that $\delta(f)$ is homogeneous of degree $m-1$. 
 
 Homogeneity condition can be weakened as follows: notice that $T\Delta^{(\CS)}_*L$ is  filtered  (as a coalgebra)
 by
 \[
 F_kT\Delta^{(\CS)}_*L=\oplus_{i\leq k}\otimes_{[i]}\Delta^{(\CS)}_*L,
 \]
Then given any $f\in Hom_{R^{\CS}}(T\Delta^{(\CS)}_*L,\Delta^{(\CS)}_*L)=\prod_{n\geq 1} Hom_{R^{[n]}}(\otimes_{[n]}\Delta^{(\CS)}_*L,\Delta^{(\CS)}_*L)$,
we similarly define $\delta(f)$ to be $\sum_{n\geq 1}\delta(f^{(n)})$, where $f^{(n)}$ is the restriction of $f$ to $\otimes_{[n]}\Delta^{(\CS)}_*L$.
This is clearly  an endomorphism of $T\Delta^{(\CS)}_*L$ that preserves the filtration.

Denote by $Coder_{\leq}(T\Delta^{(\CS)}_*L)$ the space of filtration preserving coderivations of $T\Delta^{(\CS)}_*L$.
\begin{lem}
\label{class-coderi-chi-tenso}
The map
\[
Hom_{R^{\CS}}(T\Delta^{(\CS)}_*L,\Delta^{(\CS)}_*L)\longrightarrow Coder_{\leq}(T\Delta^{(\CS)}_*L),\; f\mapsto \delta(f)
\]
is an isomorphism.
\end{lem}
With all the technology in place, the proof is no different from the ordinary one. Any filtration preserving coderivation $\delta$ is the sum
$\delta_0+\delta_1+\delta_2+\cdots$, where $\delta_m$ is a  homogeneous coderivation of degree $m-1$, and then an obvious inductive argument
shows that each $\delta_m$ is determined by its restriction to $\otimes_{[m]}L$ via formula (\ref{const-code-tenso-chira}). For example,
\[
(\delta_m\otimes 1+1\otimes\delta_m)\circ\Delta^{(\CS)}(\Delta)(x_1\cdots x_{m+1}\otimes P)=
\delta_m(x_1\cdots x_m)\otimes x_{m+1}\otimes P+x_1\otimes\delta_m(x_2\cdots x_{m+1})\otimes P,
\]
which forces 
\[
\delta_m(x_1\cdots x_{m+1}\otimes P)=\delta_m(x_1\cdots x_m)x_{m+1}\otimes P+x_1\delta_m(x_2\cdots x_{m+1})\otimes P,
\]
etc. $\qed$

\subsubsection{ }
\label{free-coal-symme-chir}
Now to the symmetric algebra case. 
We have seen already that the spaces of operations, such as $P^*_I(\{L\},L)$, carry an action of the permutation group, sect.~\ref{tansl-invar-version}.  This, of course,
has a version in $Mod(R^{\CS})$. Namely, for each $\sigma\in S_n$, there is a morphism
\[
\sigma: \otimes_{[n]}\Delta^{(\CS)}_*L\longrightarrow \otimes_{[n]} \Delta^{(\CS)}_*L
\]
determined, due to (\ref{exte-morphi-3}), by the composite morphism
\[
\otimes_{[n]}L\longrightarrow \Delta^{(\sigma)}_*\otimes_{[n]}L\hookrightarrow  (\Delta^{(\CS)}_*L)_{([n])},\;
x_1x_2\cdots x_n\otimes P\mapsto x_{\sigma_1}x_{\sigma_2}\cdots x_{\sigma_n}\otimes^\sigma P.
\]
The superscript $\sigma$ is designed simply to emphasize that  $x_{\sigma_1}x_{\sigma_2}\cdots x_{\sigma_n}\otimes^\sigma P$ is an element
of $\Delta^{(\sigma)}_*\otimes_{[n]}L$, and not of, say, $\otimes_{[n]}L$. 

If $L$ is graded, then this action can -- and will --be replaced with
 \[
\sigma: \otimes_{[n]}\Delta^{(\CS)}_*L\longrightarrow \otimes_{[n]} \Delta^{(\CS)}_*L, 
x_1x_2\cdots x_n\otimes P\mapsto \epsilon(\sigma,\vec{x})x_{\sigma_1}x_{\sigma_2}\cdots x_{\sigma_n}\otimes^\sigma P;
\]
here and elsewhere we freely use the notation introduced in sects.~\ref{lie-infty}, \ref{def-via-coderiv}.

Define the symmetric algebra $S^\bullet\Delta^{(\CS)}_*L$ to be either the quotient of $T\Delta^{(\CS)}_*L$ by the 2-sided ideal generated by
elements $a\otimes P-\sigma(a\otimes P)\in (\otimes_{[n]}\Delta^{(\CS)}_*L)_{([n])}$ or, equivalently, to be the image of the morphism
$\oplus_{n\geq 1}\sum_{\sigma\in S_n}\sigma$.  As in sect.~\ref{star-comulti-tenso-case},  $S^\bullet\Delta^{(\CS)}_*L$ is a coalgebra (in $Mod(R^\CS)$),
the coproduct being  $\Delta^{(\CS)}_*(\Delta)$, where $\Delta: S^\bullet L\rightarrow S^\bullet\otimes S^\bullet L$ is the ordinary coproduct defined in
(\ref{def-coalg-usu-sy-alg}). An explicit formula is unsurprisingly similar to the ordinary one:
\begin{equation}
\label{def-coalg-star-unusu-sy-alg}
\Delta^{(\CS)}_*\Delta(x_1x_2\cdots x_n\otimes P)=\sum_i\sum_\sigma \epsilon(\sigma,\vec{x})x_{\sigma_1}x_{\sigma_2}\cdots x_{\sigma_i}\otimes x_{\sigma_{i+1}}\cdots x_{\sigma_n}\otimes^\sigma P,
\end{equation}
the sum being extended to all $(i,n-i)$-unshuffles $\sigma$. The simplest way to come to grips with this formula is to notice that identifying $S^\bullet\Delta^{(\CS)}_*L$ with a
subobject of $T\Delta^{(\CS)}_*L$,  the coproduct just defined coincides with the restriction of the one defined in sect.~\ref{star-comulti-tenso-case}. Therefore, we have an embedding
of coalgebras $S^\bullet\Delta^{(\CS)}_*L\hookrightarrow T\Delta^{(\CS)}_*L$.

The classification of coderivations is analogous to that in the tensor algebra case.  Denote by $SP^*_{[n]}(\{L\},L)\subset P^*_{[n]}(\{L\},L)$ the subspace
of symmetric, i.e., fixed under the action of $S_n$, operations. One has an isomorphism
\[
Hom_{R^\CS}(S^\bullet\Delta^{(\CS)}_*L,L)\iso\prod_{n=1}^\infty SP^*_{[n]}(\{L\},L).
\]
\begin{sloppypar}
This shows that  $Hom_{R^\CS}(S^\bullet\Delta^{(\CS)}_*L,L)\subset Hom_{R^\CS}(T^\bullet\Delta^{(\CS)}_*L,L)$; therefore, given $f\in Hom_{R^\CS}(S^\bullet\Delta^{(\CS)}_*L,L)$,
we can consider $\delta(f)$, which is a coderivation of  $T^\bullet\Delta^{(\CS)}_*L$ by Lemma~\ref{class-coderi-chi-tenso}, hence a coderivation of $S^\bullet\Delta^{(\CS)}_*L$.
A moment's thought will show that , similarly to  (\ref{form-fo-code-ordina-ca}), if $f\in SP^*_{[n]}(\{L\},L)$, then
\end{sloppypar}
\begin{equation}
\label{form-fo-code-sta-operat-ca}
f(x_1\cdots x_N\otimes P)=\sum_\sigma \epsilon(\sigma,\vec{x})f(x_{\sigma_1},\cdots x_{\sigma_n})x_{\sigma_{n+1}},...,x_{\sigma_N}\otimes^\sigma P,
\end{equation}
the sum being extended to all $(i,n-i)$ unshuffles.

Denote by $Coder_{\leq}(S^\bullet\Delta^{(\CS)}_*L)$ the space of filtration preserving coderivations of $S\Delta^{(\CS)}_*L$.

\begin{lem}
\label{class-coderi-chi-symmetri}
The map
\[
Hom_{R^{\CS}}(S^\bullet\Delta^{(\CS)}_*L,\Delta^{(\CS)}_*L)\longrightarrow Coder_{\leq}(S^\bullet\Delta^{(\CS)}_*L),\; f\mapsto \delta(f)
\]
is an isomorphism.
\end{lem}
\subsection{ }
\label{simi-class-coalge-endo-bd}
Arguing similarly, one shows that the space of filtration preserving coalgebra morphisms  $S^\bullet\Delta^{(\CS)}_*L\rightarrow S^\bullet\Delta^{(\CS)}_*M$ is also identified with
$Hom_{R^{\CS}}(S^\bullet\Delta^{(\CS)}_*L,\Delta^{(\CS)}_*M)$. An endomorphism associated with $f=\{f_n\}$, $f_n\in SP^*_{[n]}(\{L\},L)$, is defined by the following
twin of (\ref{expl-descr-morph-infty}):

\begin{multline}
\label{expl-descr-morph-infty-in-bd-ca}
f(x_1\cdots x_n\otimes P)=\\
\sum_{1\leq i_1<i_2<i_3<\cdots <n}\sum_\sigma \epsilon(\sigma,\vec{x})f_{i_1}(x_{\sigma_1},...,x_{\sigma_{i_1}})f_{i_2-i_1}(x_{\sigma_{i_1+1}},...,x_{\sigma_{i_2}})
\cdots f_{n-i_{k}}(x_{\sigma_{i_{k}+1}},...,x_{\sigma_{n}})\otimes^\sigma P,
\end{multline}
\subsubsection{ }
\label{using-RS-symm-alg-to-descri-star-infty-li}
The relation of the just described coalgebra approach to $\text{Lie}^*_\infty$ algebras is now easy to describe along the lines of sect.~\ref{def-via-coderiv}.

Given $f\in P^*_{[n]}(\{L\},L)$ define $\hat{f}=(-1)^{n(n-1)/2}s^{-1}\circ f\circ s^{\otimes n}\in P^*_{[n]}(\{L[1]\},L[1])$. If $f$ is antisymmetric, then $\hat{f}$ is symmetric,
hence defines a coderivation $\delta(\hat{f})$ of $S^\bullet \Delta^{(\CS)}_*L[1]$.

Any $\text{Lie}^*_\infty$ algebra $(L,\{l_n\})$ defines, therefore, a coderivation, $\sum_{n\geq 1}\delta(\hat{f}_n)$, of $S^\bullet \Delta^{(\CS)}_*L[1]$. 
\begin{lem}
\label{liestarinfty-vs-coderi}
This construction sets up a 1-1 correspondence between $\text{Lie}^*_\infty$ algebra structures on $L$ and coderivations of  $S^\bullet \Delta^{(\CS)}_*L[1]$ of 
degree 1 and square 0.
\end{lem}

The proof is  identical with the proof of the corresponding result in the ordinary case, sect.~\ref{def-via-coderiv}, and will be omitted.

\subsubsection{ }
\label{defni-morph-listarinfty-case}This result will be essential for us  in that it prompts the following {\em definition}:

A $\text{Lie}^*_\infty$ morphism $(L,\{l_n\})\rightarrow (L',\{l'_n\})$  is a  morphism  $(S^\bullet \Delta^{(\CS)}_*L[1],\delta), (S^\bullet \Delta^{(\CS)}_*L'[1],\delta')$
of coalgebras with coderivation.

In other words, it is an $f$ defined by (\ref{expl-descr-morph-infty-in-bd-ca}) that  satisfies $\delta'\circ f=f\circ\delta$.

\subsection{ }
\label{the-real-thing}
The discussion above is but a shadow of the genuine  Beilinson-Drinfeld category \cite{BD},  2.2.  Given a smooth algebraic  curve $X$, their category is one of
right $D_X$-modules with the pseudo-tensor structure defined by 
\[
P^*_I(\{\CM_i\},\CN)=\text{Hom}_{D_{X^I}}(\boxtimes_I \CM_i\longrightarrow \Delta_\ast \CN),
\]
where $\Delta: X\longrightarrow X^I$ is the diagonal embedding.

Seeking to spell out everything in the simplest possible case, let from now on $X$ be $\BC$, $X^I=\times_I X$, $\BC[X^I]$ the corresponding polynomial ring, $D_{X^I}$
the corresponding algebra of globally defined differential operators; we let $x$ be the coordinate on $X$, $\partial_x=\partial/\partial x$. The  various products
over an arbitrary finite set, here and elsewhere,  are made sense of along the lines of   sect.~\ref{from-n-to-any-set}

Given a surjection $\pi: J\twoheadrightarrow I$, there arise an embedding $X^I \hookrightarrow X^J$ and  the corresponding algebra homomorphism  $\BC[X^J]\twoheadrightarrow
\BC[X^I]$, $x_j\mapsto x_{\pi(j)}$. Define $D_{X^I\rightarrow X^J}=\BC[X^I]\otimes_{\BC[X^J]}D_{X^J}$, which is operated on by $D_{X^J}$ on the right -- obviously, and
by $D_{X^I}$ on the left via $\partial_{x_i}\mapsto\sum_{j\in\pi^{-1}(i)}\partial_{x_j}$; this makes $D_{X^I\rightarrow X^J}$ into a $D_{X^I}-D_{X^J}$-bimodule. There are obvious
isomorphisms: 
\[
\otimes_I D_X\iso D_{X^I}, \; \otimes_I D_{X\rightarrow X^{J_i}}\iso D_{X^I\rightarrow X^J}, \;
D_{X^I\rightarrow X^J}\otimes_{D_{X^J}}D_{X^J\rightarrow X^K}\iso D_{X^I\rightarrow X^J};
\]
$J_i$ stands for $\pi^{-1}(i)$ in the 2nd isomorphism.

For a collection of right $D_X$-modules, $\CM_i$, $i\in I$, $\CN$, define
\[
P^*_I(\{\CM_i\},\CN)=\text{Hom}_{D_{X^I}}(\otimes_I \CM_i\longrightarrow \CN\otimes_{D_X}D_{X\rightarrow X^I}),
\]
The composition is defined as follows:  for a surjection $\pi: J\twoheadrightarrow I$, and a collection of operations
$\psi_i\in P^*_{J_i}(\{\CL_j\}, \CM_i)$, $i\in I$, $J_i=\pi^{-1}(i)$, and $\phi\in P^*_I(\{\CM_i\},\CN)$, define $\phi(\psi_i)\in P^*_J(\{\CL_j\},\CN)$ to be the 
composite map, cf. sect.~\ref{tansl-invar-version}:
\begin{multline}
\otimes_J\CL_j=\otimes_I\otimes_{J_i}\CL_j\stackrel{\otimes\psi_i}{\longrightarrow}\otimes_I (\CM_i\otimes_{D_X} D_{X\rightarrow X^{J_i}})=(\otimes_I \CM_i)\otimes_
{D_{X^I}} D_{X^I\rightarrow X^{J}}\\
\stackrel{\phi}{\longrightarrow} (\CN\otimes_{D_X}D_{X\rightarrow X^I})\otimes_{D_{X^I}} D_{X^I\rightarrow X^{J}}
\iso \CN\otimes_{D_X}D_{X\rightarrow X^J}.\nonumber
\end{multline}
The associativity follows from the isomorphisms $D_{X^I\rightarrow X^J}\otimes_{D_{X^J}}D_{X^J\rightarrow X^K}\iso D_{X^I\rightarrow X^J}$.

\subsection{ }
\label{cm-vs-cmd-some-disc}
Denote by $\CM_D^*$ the pseudo-tensor category just defined.  Just as $\CM^*$ of sect.~\ref{tansl-invar-version}, in fact as any pseudo-tensor category, it carries commutative associative, Lie, Poisson,
etc. objects, which we will still be calling commutative$^!$, Lie, coisson, etc., algebras.

The obvious similarity between $\CM_D^*$ and $\CM^*$
is easily made into an assertion as follows. Given an $R$-module $M$, $M[x]$ is naturally a $D_X$-module if we stipulate $m\partial_x=m\partial$, $m\in M$.
This defines a functor
\[
\Phi:\CM^*\longrightarrow \CM_D^*,\; M\mapsto M[x],
\]
which is clearly pseudo-tensor and faithful. In fact,
it  identifies $\CM^*$ with the {\em translation-invariant} subcategory of
 $\CM^*_D$, i.e., $\CL$ is isomorphic to $\Phi(M)$ for $M\in\CM^*$ precisely when $\CL$ is translation-invariant, and  $\phi\in P^*_I(\{M_i[x]\},N[x])$
 belongs to $\Phi P^*_I(\{M_i\},N)$ if and only if $\phi$ is translation-invariant. Therefore, an object of some type of $\CM^*$ is the same as a translation-invariant object
 of the same type in $\CM^*_D$. 
 
 \subsection{ }
 \label{Lie-alf-fourier-comp}
 We are exclusively interested in the translation invariant objects, but even then this more general point of view is helpful. The assignment $\CM^*_D\ni\CL\mapsto
 h(\CL)\stackrel{\text{def}}{=}\CL/\CL\partial_x\in\CV ect$ is still an augmentation functor, sect.~\ref{augm-fnctr}, and if $L$ is a Lie* algebra, then
  $h(L[x])$ is a Lie algebra, just as $h(L)$, 
 sect.~\ref{from lie*tolie}.
 
 Likewise, since our discussion easily localizes, if $L$ is a Lie* algebra, then $h(L[x,x^{-1}])$ is a Lie algebra. If we let $a_{[n]}$ denote the class of
 $a\otimes x^n$ in $h(L[x,x^{-1}])$, then it is immediate to derive from  (\ref{jacobi}) a formula for the bracket:
\begin{equation}
\label{borch-in-hl} 
[a_{[n]},b_{[m]}]=\sum_{j\geq 0}{n\choose j}(a_{(j)}b)_{[n+m-j]},
\end{equation}
also called the Borcherds identity. Denote this Lie algebra $\text{Lie}(L)$; cf. \cite{K}, pp.41-42, \cite{FBZ}, 16.1.16.

\subsection{ }
\label{defn-chir-alg}
Similarly, the concept of a chiral algebra, even in the translation-invariant setting, is most naturally introduced in the framework of
$D_X$-modules.  For an $I$-family $\{\CA_i\}\subset \CM_D$ let  $\Delta_{\alpha\beta}=(x_\alpha-x_\beta)^{-1}$, $\alpha,\beta\in I$ 
and denote  by  $\otimes_I\CA_i[\{ \Delta_{\alpha\beta}\}]\in \CM_D$
the localization at the indicated elements. Define
\[
P^{ch}_I(\{\CA_i\},\CN)=\text{Hom}_{D_{X^I}}(\otimes_I\CA_i[\{ \Delta_{\alpha\beta}\}],\CN\otimes_{D_X}D_{X\rightarrow X^I}).
\]
Elements of such sets are called {\em chiral} operations. They are composed in the same way as the *-operations of sect.~\ref{the-real-thing}, except that
now one has to deal with the poles. Let us examine the simplest and most important
such composition; the pattern will then become clear.

Fix $\phi\in\text{Hom}_{D_{X^2}}(\CL\otimes\CL[\Delta], \CL)$. The composition $\phi(\phi,id)$ is defined as follows:
\begin{multline}
\CL\otimes\CL\otimes\CL[\Delta_{12},\Delta_{13},\Delta_{23}]\iso(\CL\otimes\CL[\Delta_{12}])\otimes\CL[\Delta_{13},\Delta_{23}]\stackrel{\phi\otimes id}{\longrightarrow}
(\CL\otimes_{D_X}D_{X\rightarrow X^2})\otimes\CL[\Delta_{13},\Delta_{23}]\\
\iso(\CL\otimes\CL)\otimes_{D_{X^2}}D_{X^2\rightarrow X^3}[\Delta_{13},\Delta_{23}]\iso(\CL\otimes\CL[\Delta])\otimes_{D_{X^2}}D_{X^2\rightarrow X^3}\stackrel{\phi\otimes id}{\longrightarrow}\\
(\CL\otimes_{D_X}D_{X\rightarrow X^2})\otimes_{D_{X^2}}
D_{X^2\rightarrow X^3}\iso \CL\otimes_{D_X}D_{X\rightarrow X^3}\nonumber.
\end{multline}
In this composition, the middle isomorphism in the second line follows from the fact that $D_{X^2\rightarrow X^3}$ is a $\BC[X^2]-\BC[X^3]$-bimodule
and by definition 
\[
D_{X^2\rightarrow X^3}[\Delta_{13},\Delta_{23}]\iso [\Delta]D_{X^2\rightarrow X^3}.
\]

This gives the category of right $D_X$-modules another pseudo-tensor structure, to be denoted $\CM^{ch}_D$. 

\subsection{ }
\label{funct-real-dir-imag}
Note a useful isomorphism of right $D_{X^2}$-modules
\[
\CL\otimes_{D_X}D_{X\rightarrow X^2}\iso \Omega^1_X\otimes\CL[\Delta]/\Omega^1_X\otimes\CL,\;
l\otimes 1\mapsto \frac{dx\otimes l}{x_1-x_2} \text{ mod } \Omega^1_X\otimes\CL,
\]
which is a manifestation of the Kashiwara lemma, \cite{Bor}, 7.1. Notice that from this point of view, the composite map
\[
\CL\otimes_{D_X}D_{X\rightarrow X^2}\twoheadrightarrow\CL\otimes_{D_X}D_{X\rightarrow X^2}/(\CL\otimes_{D_X}D_{X\rightarrow X^2})\partial_1\iso \CL
\]
is defined by the residue
\[
\Omega^1_X\otimes\CL[\Delta]\longrightarrow\CL,\; \omega\otimes l\mapsto (\int_{x_1:|x_1-x_2|=r}\omega)l.
\]

\subsection{ }
\label{main-defn-chir-case}
A Lie$^{ch}$ algebra (on $X$) is a Lie object in $\CM^{ch}_D$; explicitly, it is a right $D_X$-module $\CL$ with  chiral bracket $[.,.]^{ch}\in P^{ch}_2(\{\CL,\CL\},\CL)$ that is 
anticommutative and satisfies the Jacobi identity.  The simplest example is $\Omega^1_X$ with the canonical right $D_X$-module structure (given by the negative
Lie derivative) and the chiral Lie bracket
\begin{equation}
\label{omega-as-chir-alg}
\Omega^1_X\otimes\Omega^1_X[\Delta]\twoheadrightarrow\Omega^2_{X^2}[\Delta]/\Omega^2_{X^2}
\iso\Omega^1_X\otimes_{D_X}D_{X\rightarrow X^2},\nonumber
\end{equation}
where the rightmost isomorphism has just been discussed, sect.~\ref{funct-real-dir-imag}. Note that the anticommutativity follows from
the fact that the natural $S_2$-equivariant structures of $\Omega^1_X\otimes\Omega^1_X$ and $\Omega^2_{X^2}$ differ
by the sign representation of $S_2$.

The chiral algebra   is  a Lie$^{ch}$ algebra  $\CL$ with  a unit, i.e., a morphism $\iota:\Omega^1_X\longrightarrow\CL$ s.t. the composition $[\iota(.),.]$
coincides with the map
\begin{equation}
\label{another-manif-kashiw}
\Omega^1_X\otimes\CL[\Delta]\twoheadrightarrow\Omega^1_X\otimes\CL[\Delta]/\Omega^1_X\otimes\CL
\iso\CL\otimes_{D_X}D_{X\rightarrow X^2},\nonumber
\end{equation}

The obvious map $\otimes_I\CA_i\longrightarrow \otimes_I\CA_i[\cup \Delta_{\alpha\beta}]$ defines, by restriction, a map
\begin{equation}
\label{how-chira-ope-def-star-opera}
P^{ch}_I(\{\CA_i\},\CN)\longrightarrow P^*_I(\{\CA_i\},\CN),
\end{equation}
hence a {\em forgetful }functor $\CM^{ch}_D\longrightarrow \CM^*_D$. It follows that each Lie$^{ch}$ algebra can be regarded as a Lie* algebra. Further composing with
$h:\CM^*_D\rightarrow \CV ect$, sect.~\ref{Lie-alf-fourier-comp}, will attach an ordinary Lie algebra $h(\CL)$ to each chiral algebra $\CL$.

A chiral algebra is called {\em commutative}
if the corresponding Lie* algebra is abelian, i.e., the corresponding Lie* bracket is 0. In the translation-invariant setting, a commutative chiral algebra is the same thing as an
ordinary unital commutative associative algebra with derivation; we shall have more to say on this in sect.~\ref{expl-formulas-borch-etc}.

The definition of a chiral algebra module  should be evident; any chiral algebra module is automatically a module over the corresponding Lie* algebra.
If $\CL$ is a chiral algebra and $\CM$ an $\CL$-module, then $h(\CL)$ is a Lie algebra,
and both $\CM$ and $h(\CM)$ are $h(\CL)$-modules. If the structure involved is translation invariant, in particular,
$\CL=L[x]$, $\CM=M[x]$, then  the fiber  $M$ is also an $h(\CL)$-module, as well as $\text{Lie}(\CL)$-module, see sect.~\ref{Lie-alf-fourier-comp}.

\subsection{ }
\label{comm-chir-mod-vs-comm-der-mod}
$\CM$, a module over a chiral algebra $\CL$, is called {\em central} if it is trivial over the corresponding Lie* algebra, \cite{BD}, 3.3.7.

In view of what is said at the end of sect.~\ref{main-defn-chir-case} it may sound as a surprise that a module over a commutative chiral algebra $\CL=L[x]$  is
not the same thing as a module over $L$ regarded as a commutative associative algebra with derivation. However, if the module in question is central, then
the two notions coincide; we shall explain this in sect.~\ref{expl-formulas-borch-etc} and show an example in sect.~\ref{ex-poly-ring}.

\subsection{ }
\label{ope}
An explicit description of a chiral algebra usually arises in the following situation. Let $\CV$ be a translation-invariant {\em left} $D_X$-module, which
amounts to having $\CV=V[x]$, $V$ being a {\em left} $R$-module. Let $\CV^r\stackrel{\text{def}}{=}\CV\otimes_{\BC[x]}\Omega^1_X$ be the corresponding
{\em right} $D_X$-module; we shall sometimes write simply $V[x]dx$ for $\CV^r$.

Notice canonical isomorphisms of right $D_{X^2}$-modules
\begin{multline}
\label{expl-descr-dir-im-codim1}
\CV^r\otimes_{D_X}D_{X\rightarrow X^2}\longrightarrow \BC[x]\otimes V[y][(x-y)^{-1}]dx\wedge dy/\BC[x]\otimes V[y]dx\wedge dy\\
\longleftarrow V[x]\otimes \BC[y][(x-y)^{-1}]dx\wedge dy/V[x]\otimes \BC[y]dx\wedge dy;
\end{multline}
the first is discussed in sect.~\ref{funct-real-dir-imag}, the second is the result of a formal Taylor series expansion
\[
v(x)\mapsto \sum_{n=0}^{\infty}\frac{1}{n!} \partial^nv(y)(x-y)^n,
\]
which is essentially Grothendieck's definition of a connection.

In this setting, the translation-invariant chiral bracket $[.,.]\in P^{ch}_2(\{\CV^r,\CV^r\},\CV^r)$ is conveniently encoded by a map, usually  referred to as an OPE:
\begin{equation}
\label{form-for-ope}
V\otimes V\longrightarrow V((x-y)),\; a\otimes b\mapsto \sum_{n\in\BZ}a_{(n)}b\otimes (x-y)^{-n-1},\; a_{(n)}b\in V.
\end{equation}
Given an OPE, one recovers the chiral bracket
\[
(V[x]dx\otimes V[y]dy)[(x-y)^{-1}]\longrightarrow \BC[x]\otimes V[y][(x-y)^{-1}]dx\wedge dy/\BC[x]\otimes V[y]dx\wedge dy
\]
by defining
\[
\frac{1}{(x-y)^N}(a\otimes f(x))\otimes(b\otimes g(y))\mapsto \sum_{n\in\BZ} a_{(n)}b\frac{f(x)g(y)}{(x-y)^{N+n+1}}dx\wedge dy \text{ mod reg},
\]
``mod reg.'' meaning, of course, ``modulo $\BC[x]\otimes V[y]dx\wedge dy$.''
In fact, this sets up a 1-1 correspondence between binary chiral operations and OPEs, \cite{FBZ}, 19.2.11, or \cite{BD}, 3.5.10.

In this vein, the Jacobi identity can also be made explicit.  The main diagonal in $X^3$ being of codimension 2,  $\CV^r\otimes_{D_X}D_{X\rightarrow X^3}$ does not allow a description
as simple as (\ref{expl-descr-dir-im-codim1}), and one relies instead on iterations of (\ref{expl-descr-dir-im-codim1}).  Writing $\CV^r\otimes_{D_X}D_{X\rightarrow X^3}$
as $(\CV^r\otimes_{D_X}D_{X\rightarrow X^2})\otimes_{D_{X^2}}D_{X^2\rightarrow X^3}$, which requires a choice of an embedding $X^2\rightarrow X^3$, such as
$(u,v)\mapsto (u,v,u)$, one obtains identifications, such as
\[
\CV^r\otimes_{D_X}D_{X\rightarrow X^3}\longrightarrow \BC[x]\otimes(\BC[y]\otimes V[z][(y-z)^{-1}]\text{ mod reg. })[(x-z)^{-1}] \text{ mod reg.};
\]
we omit differentials, $dx\wedge dy\wedge dz$, for typographical reasons.

Write  $a(x-y)b$ for OPE (\ref{form-for-ope}).  Various compositions that enter the Jacobi identity involve expressions such as
\begin{eqnarray}
&&[a,[b,c]]=(a(x-z)(b(y-z)c \text{ mod reg.})\text{ mod reg. })dx\wedge dy\wedge dz, \label{wh-reg-wh-arent-1}\\
&&[[a,b]c]=( (a(x-y)b \text{ mod reg.})(y-z)c\text{ mod reg. })\,dx\wedge dy\wedge dz, \text{ etc.}.\nonumber
\end{eqnarray}
The Jacobi identity,
\[
[a,[b,c]]-[a,[b,c]]^{(1,2)}-[[a,b],c]=0,
\]
implies that for any $F(x,y,z)=(x-y)^r(x-z)^s(y-z)^t$, $r,s,t\in\BZ$,
\begin{multline}
\label{borcherds-full-ident}
\oint_{x:|x-z|=R}dx\oint_{y:|y-z|=r}dy F(x,y,z)a(x-z)b(y-z)c\,dz\\
-\oint_{y:|y-z|=R}dy\oint_{x:|x-z|=r}dx F(x,y,z)b(y-z)a(x-z)c\,dz\\
-\oint_{y:|y-z|=R}dy\oint_{x:|x-y|=r}dx F(x,y,z)(a(x-y)b)(y-z)c\,dz=0,
\end{multline}
where $R>r$.  Let us explain this.

Denote by
$Jac\in P^{ch}_3(\{\CV^r,\CV^r,\CV^r\},\CV^r)$ the left hand side of the Jacobi identity; it is a map

\begin{equation}
Jac:( V[x]\otimes V[y]\otimes V[z])[(x-y)^{-1},(x-z)^{-1},(y-z)^{-1}]dx\wedge dy\wedge dz\longrightarrow V[t]dt\otimes_{D_X}D_{X\rightarrow X^3}.\nonumber
\end{equation}

  Written down in terms of the OPE it gives the left hand side of (\ref{borcherds-full-ident}) except that:  
  
  the $\int$ signs must be removed;
  
 the function $F(x,y,z)$ must be  expanded in powers of appropriate variabes, $(x-z)$ and $(y-z)$ for the 1st and 3rd term, $(x-y)$ and $(y-z)$ for the 2nd one, in domains
prescribed by  the definition of the composition, sect.~\ref{defn-chir-alg}; for example, in the case of the 1st integral, one has
\[
F(x,y,z)=(x-y)^r(x-z)^s(y-z)^t=(x-z)^{s+r}(y-z)^t\sum_{j=0}^{\infty} (-1)^j {r\choose j}\left(\frac{y-z}{x-z}\right)^j;
\]

finally, regular terms must be crossed out, see (\ref{wh-reg-wh-arent-1}). 

Treating the arising 3 expressions requires an effort as they belong to 3 different realizations of the same space,
$V[t]dt\otimes_{D_X}D_{X\rightarrow X^3}$. However, part of this computation is easy: the composition
\begin{multline}
\nonumber
( V[x]\otimes V[y]\otimes V[z])[(x-y)^{-1},(x-z)^{-1},(y-z)^{-1}]dx\wedge dy\wedge dz\stackrel{Jac}{\longrightarrow} V[t]dt\otimes_{D_X}D_{X\rightarrow X^3}\\
\twoheadrightarrow V[t]dt =V[t]dt\otimes_{D_X}D_{X\rightarrow X^3}/(V[t]dt\otimes_{D_X}D_{X\rightarrow X^3})\BC[\partial_x,\partial_y]
\end{multline}
is defined simply by taking the residues, just as in  sect.~\ref{funct-real-dir-imag}, hence it equals the left hand side of (\ref{borcherds-full-ident}).

Formula (\ref{borcherds-full-ident}) is the Borcherds identity \cite{Borch} in the form suggested in \cite{K}, 4.8. Therefore, a translation invariant chiral algebra on $\BC$ defines a vertex algebra. A passage
in the opposite direction is carefully explained in \cite{FBZ}, Ch.15. Here is an alternative argument:
  the image of $Jac$ is a $D_{X^3}$-submodule, and if (\ref{borcherds-full-ident}) is valid, then this submodule
  belongs to $Im\partial_x+ Im\partial_y$, hence equals 0 according to the Kashiwara lemma, as desired.
  
  Originally, the comparative analysis of the notions of chiral and vertex algebra was carried out in \cite{HL} .

\subsection{ }
\label{expl-formulas-borch-etc}
The case $F(x,y,z)=(x-z)^m(y-z)^n$ of (\ref{borcherds-full-ident}) reproduces the Borcherds commutator formula (\ref{jacobi})
\begin{equation}
\label{borch-comm-again}
 a_{(n)}b_{(m)}c-b_{(m)}a_{(n)}c=\sum_{j\geq 0}{n\choose j}(a_{(j)}b)_{(n+m-j)}c.
 \end{equation}
 The case $F(x,y,z)=(x-y)^{-1}(y-z)^{n}$ becomes the celebrated normal ordering formula
 \begin{equation}
 \label{norm-order-form}
 (a_{(-1)}b)_{(n)}c=\sum_{j=0}^\infty b_{(-j+n-1)}a_{(j)}c+a_{(-1)}b_{(n)}c+\sum_{j=0}^\infty a_{(-j-2)}b_{(j+n+1)}c
 \end{equation}
In fact, these particular cases suffice to reproduce the entire (\ref{borcherds-full-ident}),  \cite{K}, 4.8.

One sees at once that in this language  a (translation-invariant) chiral algebra (on $\BC$) is a vector space $V$ with a family of multiplications, $_{(n)}$, $n\in\BZ$.   The unit axiom (\ref{another-manif-kashiw}) reads: there is an element $\Vac$ such that $a(z)\Vac\equiv a\text{ mod }(z)$
and $\partial\Vac=0$.

``$V$ is commutative''
(see sect.~\ref{main-defn-chir-case}) means the ``$n$th product is 0 if $n\geq 0$.''  If so,
 (\ref{norm-order-form}) with $n=-1$ shows that the product $_{(-1)}$ is associative, and then
Borcherds commutator formula (\ref{borch-comm-again}) shows that it is commutative, and so $(V,_{(-1)},\Vac,\partial)$ is an associative, commutative, unital algebra with derivation. The passage in the opposite direction is explained in \cite{K,FBZ}.

Similarly, the conceptual definition of a chiral algebra module, reviewed in sect.~\ref{main-defn-chir-case}, boils down to a vector space $M$ with multiplications
\[
_{(n)}^M:\, V\otimes M\longrightarrow M,\; a\otimes m\mapsto a^M_{(n)}m
\]
so that
\begin{eqnarray}
 & &a_{(n)}^Mb_{(m)}^Mc-b_{(m)}^Ma_{(n)}^Mm=\sum_{j\geq 0}{n\choose j}(a_{(j)}b)_{(n+m-j)}^Mm,\\
 \label{norm-order-form-module}
 & &(a_{(-1)}b)_{(n)}^Mc=\nonumber\\
 & &\sum_{j=0}^\infty b_{(-j+n-1)}^Ma_{(j)}^Mc+a_{(-1)}^Mb_{(n)}^Mc+\sum_{j=0}^\infty a_{(-j-2)}^Mb_{(j+n+1)}^Mm;
 \end{eqnarray}
 we are deliberately omitting some of the obvious axioms.
 
 One sees clearly how the concept of a module over a commutative chiral algebra is different from one over a commutative associative algebra with derivation:
 the associativity condition $(ab)m=a(bm)$ in the latter is replaced by the more cumbersome  (\ref{norm-order-form-module}) in the former. If, however, $M$ is central, sect.~\ref{comm-chir-mod-vs-comm-der-mod},  which means that $a^M_{(n)}m=0$, $n\geq 0$, then the ``correction terms'' in (\ref{norm-order-form-module}) vanish, and the two concepts
 become indistinguishable. 
 
 \subsection{ }
 \label{vert-envel-lie-star}
 We have seen, sect.~\ref{main-defn-chir-case}, that there is a forgetful functor that makes a chiral algebra into a Lie* algebra. This functor admits the left adjoint called
 the {\em chiral enveloping} algebra. Let us sketch its construction, cf. \cite{FBZ}, 16.1.11, \cite{BD}, 3.7.1. (We work in the translation-invariant setting, this goes without saying.)
 
 Given a Lie* algebra $L$, consider the Lie algebra $\text{Lie}(L)$, sect.~\ref{Lie-alf-fourier-comp}. Formula (\ref{borch-in-hl}) implies that  $\text{Lie}(L)_+$ defined to
 be spanned by $a_{[n]}$, $a\in L$, $n\geq 0$, is a Lie subalgebra. Define $U^{ch}L$ to be $U(\text{Lie}(L))/U(\text{Lie}(L)_+)$.  Here $U(.)$ is the ordinary universal enveloping
 of a Lie algebra.
 
 It is easy to see that the  map $L\longrightarrow \text{Lie}(L)$, $a\mapsto a_{[-1]}$, is injective, and so is the composition
 \[
 L\longrightarrow \text{Lie}(L)\hookrightarrow U(\text{Lie}(L))\twoheadrightarrow U(\text{Lie}(L))/U(\text{Lie}(L)_+)
 \]
 The Reconstruction Theorem,
 \cite{FBZ}, 2.3.11 or \cite{K}, 4.5, implies that $U^{ch}L$ carries a chiral algebra structure defined, in terms of $_{(n)}$-products, by a slightly tautological formula
 \[
 (\overline{a_{[-1]}})_{(n)} v=a_{[n]}\cdot v;
 \]
 here $\overline{a_{[-1]}}$ is the image of $a_{[-1]}$ under the above composition, 
 and $\cdot$ on the right means the action of $\text{Lie}(L)$ on $U(\text{Lie}(L))/U(\text{Lie}(L)_+)$.

\section{cdo}
\label{cdo}
\subsection{ }
We shall work exclusively in the translation-invariant situation, although much of what we are about  to say does not require this assumption,
and so we shall typically deal with fibers of the actual objects, cf. sect.~\ref{cm-vs-cmd-some-disc}, \ref{expl-formulas-borch-etc}.  Thus, for example,
the phrase `` a chiral (Lie*, etc.) algebra $V$'' means the fiber of a translation-invariant chiral (Lie*, etc.) algebra $V[x]$, and a chiral (Lie*, etc.) algebra
morphism $f:\;V\longrightarrow W$ means $f\otimes  id:\; V[x]\longrightarrow W[x]$.

\subsection{ }
\label{defftn-cdo}
Let $A$ be a commutative associative unital algebra.  A chiral algebra $\CD^{ch}_A$ is called an algebra of {\em chiral differential operators} over $A$
if it carries a filtration $F_{-1}\CD^{ch}_A=\{0\}\subset F_0\CD^{ch}_A\subset F_1\CD^{ch}_A\subset\cdots$, $\cup_n F_n\CD^{ch}_A=\CD^{ch}_A$, s.t. the graded object 
\[
\text{gr}\CD^{ch}_A=\bigoplus_{n=0}^{\infty}F_n\CD^{ch}_A/F_{n-1}\CD^{ch}_A
\]
is a coisson algebra, sect.~\ref{comm-modul-act-prod-comm}, which is isomorphic, as a coisson algebra, to $J_\infty S^\bullet_A T_A$, sect~\ref{ poiss-jet-coiss}, \ref{cm-vs-cmd-some-disc}.

By definition,  $F_0\CD^{ch}_A=J_\infty A$ and $F_1\CD^{ch}_A$ fits in the short exact sequence
\begin{equation}
\label{chiral-alg-as-appear}
0\longrightarrow J_\infty A\longrightarrow F_1\CD^{ch}_A\longrightarrow J_\infty T_A\longrightarrow 0,
\end{equation}
{\em loc. cit.}  Notice that both $F_1\CD^{ch}_A$ and $J_\infty T_A$ are Lie* algebras and chiral $J_\infty A$-modules, but while $J_\infty T_A$ is a Lie* $J_\infty A$-algebroid,
$F_1\CD^{ch}_A$ is not. This has to do with the fact that $J_\infty A$ being a commutative algebra with derivation is both a commutative$^!$ algebra, sect.~\ref{comm-modul-act-prod-comm}, and a commutative
chiral algebra, sect.~\ref{main-defn-chir-case}; in its former capacity it operates on $J_\infty T_A$, but it acts on $F_1\CD^{ch}_A$ only as a chiral algebra, sect.~\ref{comm-chir-mod-vs-comm-der-mod}.

This prompts the following definition.

\subsection{ }  
\label{defn-chir-algebr}
A chiral algebroid ($A$-algebroid)\footnote{we should have said `` a translation-invariant chiral algebroid on $\BC$ in the case of a jet-scheme''}
is a short exact sequence
\begin{equation}
\label{sh-ex-seq-def-chir-alg}
0\longrightarrow J_\infty A\stackrel{\iota}{\longrightarrow} \CL_A^{ch}\stackrel{\sigma}{\longrightarrow} J_\infty T_A\longrightarrow 0,
\end{equation}
where $\CL_A^{ch}$ is a Lie* algebra and a chiral module over $J_\infty A$, and the arrows respect all the structures.  Here is what this amounts to.

(i) $J_\infty A\stackrel{\iota}{\longrightarrow} \CL_A^{ch}$ is a morphism of chiral $J_\infty A$-modules and $\text{Lie}^*$ algebras. 

(ii)  $J_\infty A$ acts as a Lie* algebra on $\CL_A^{ch}$ in two ways; first, via an adjoint action as a subalgebra of the Lie* algebra $\CL_A^{ch}$, second, because
a chiral action of $J_\infty A$ on $\CL_A^{ch}$ induces a Lie* action, sect.~\ref{main-defn-chir-case}. We require that these two actions coincide.

(iii) $\CL_A^{ch}\stackrel{\sigma}{\longrightarrow} J_\infty T_A$ is a Lie* algebra and $J_\infty A$-module morphism. 

(iv) Items (i) and (iii) imply that $\iota(J_\infty A)\subset \CL^{ch}_A$ is an abelian Lie* ideal and, therefore, is acted upon by $J_\infty T_A$. We require that this action
be equal to the canonical action of $J_\infty T_A$ on $J_\infty A$, 
sect.~\ref{ poiss-jet-coiss}.

(v) The Lie* action of $\CL^{ch}_A$ on itself is a derivation of  the chiral action of $J_\infty A$ on $\CL_A^{ch}$. Namely, 
\begin{equation}
\label{meani-of-lch-linea-defni-cdo}
[.,\mu(.,.)]=\mu([.,.],.)+\mu(.,[.,.])^{(1,2)},
\end{equation}
where $\mu\in P^{ch}_{\{1,2\}}(\{J_\infty A, \CL^{ch}_A\},\CL^{ch}_A)$ is the chiral module structure and 
$[.,.]\in P^{*}_{\{1,2\}}(\{\CL_A^{ch},\CL_A^{ch}\},\CL^{ch}_A)$ is the Lie*-bracket.

{\em Remark.}

 Point (v) is a straightforward analogue of  (\ref{lie-alg-comp-str}). In terms
of $_{(n)}$-products it amounts to the fact that the  commutator formula,  cf. sect.~\ref{expl-formulas-borch-etc},
\begin{equation}
\label{what-l-linear-means}
a_{(n)}b_{(m)}c-b_{(m)}a_{(n)}c=\sum_{j\geq 0}{n\choose j}(a_{(j)}b)_{(n+m-j)}c,
\end{equation}
whose validity for $m,n\geq 0$ is the consequence of $\CL^{ch}_A$ being a Lie* algebra, is also valid for $m<0$ if $b\in J_\infty A$, cf. (\ref{brac-vect-fields}).

\subsection{ }
\label{ex-poly-ring}
A well-known example arises when $A=\BC[x_1,...,x_N]$. Introduce $\fg$,  a Lie algebra with generators $x_{ij}$, $\partial_{mn}$, $1\in\BC$ and relations
$[\partial_{mn}, x_{ij}]=\delta_{mi}\delta_{n,-j}$.  There is a subalgebra, $\fg_{-}$, defined to be the linear span of $x_{ij}$, $\partial_{mn}$, $j>0$, $m\geq 0$.
The induced representation $\text{Ind}_{\fg_{-}}^{\fg}\BC$, which is naturally identified with $\BC[x_{ij},\partial_{mn};\;j\leq 0,n<0]$, is well known to carry
a vertex algebra structure; it is often referred to as a ``$\beta$-$\gamma$-system.  Explicit formulas can be found in \cite{MSV}.  For example, one has
\begin{multline}
\label{defnbetagammasys}
 (\partial_{i,-1})_{(0)}(x_{j,0})=\delta_{ij},\\
   (\partial_{i,-1})_{(n+1)}(x_{j,0})= (\partial_{i,-1})_{(n)}(\partial_{j,-1})= (x_{i,0})_{(n)}(x_{j,0})=0\text{ if }n\geq 0.
\end{multline}

Fix $\BC[x_1,...,x_n]\hookrightarrow \BC[x_{ij},\partial_{mn};\;j\leq 0,n<0]$, $x_i\mapsto x_{i0}$. For any  \'etale localization  $\BC[x_1,...,x_N]\subset A$,
the space 
\[
A[x_{ij},\partial_{mn};\;j,n<0]\stackrel{\text{def}}{=}A\otimes_{\BC[x_1,...,x_N]}\BC[x_{ij},\partial_{mn};\;j\leq 0,n<0]
\]
inherits a vertex algebra structure from $\BC[x_{ij},\partial_{mn};\;j\leq 0,n<0]$.

The increasing filtration $\{F_rA[x_{ij},\partial_{mn};\;j,n<0]\}$, $r\geq 0$, is defined by counting the letters $\partial_{mn}$, $n<0$. The graded object is identified with 
$J_\infty S^\bullet_A T_A$,
and so $A[x_{ij},\partial_{mn};\;j,n<0]\}$ is a CDO, sect.~\ref{defftn-cdo}.

The space $F_1A[x_{ij},\partial_{mn};\;j,n<0]$  is a chiral algebroid. Exact sequence (\ref{sh-ex-seq-def-chir-alg}) in this case becomes
\begin{equation}
0\longrightarrow A[x_{ij};\;j<0]\longrightarrow F_1 A[x_{ij},\partial_{mn};\;j,n<0]\longrightarrow 
\oplus_m\oplus_{n<0}A[x_{ij};\;j<0]\partial_{mn}\longrightarrow 0.\nonumber
\end{equation}

It is easy to see exactly how $F_1 A[x_{ij},\partial_{mn};\;j,n<0]$ fails to be a central chiral $J_\infty A$-module and $J_\infty T_A=F_1 A[x_{ij},\partial_{mn};\;j,n<0]/J_\infty A$ does not:
suppressing extraneous indices we derive using (\ref{defnbetagammasys}, \ref{norm-order-form-module})
\[
((x_0)^2)_{(-1)}\partial_{-1}=x_0^2\partial_{-1}-2x_{-1}.
\]

\subsection{ }
\label{class-chir-algebr}
Classification of chiral algebroids is delightfully similar to that of Picard-Lie algebroids, sect.~\ref{defn-p-lie}.  

To begin with, assume that
the tangent Lie algebroid $T_A$ is a free $A$-module with basis $\{\xi_i\}$. Then it is easy to see that the chiral module structure
on $\CL^{ch}_A$ cannot be deformed. Indeed,  suppose one such structure is given.  It follows at once that any element $\xi$  of
$\CL^{ch}_A$ is uniquely written as a sum: $\xi=f_0+\{\sum_{i>0} f_{i,(-1)}\partial^{j_i}\xi_i$, $f_i\in J_\infty A$.
  Elements $f_{(n)}\xi$, $n\geq0$,
are then completely determined: they form the Lie* action of $J_\infty A$ on $\CL^{ch}_A$, which by (iii) is the same as  minus the adjoint action of $\CL^{ch}_A$
restricted to $J_\infty A$, which by (ii) is the pull-back via $\CL_A^{ch}\longrightarrow J_\infty T_A$ of the canonical action of $J_\infty T_A$ on $J_\infty A$.

Finally, we have to compute elements of the type $g_{(-1)}\xi$, say, $g_{(-1)}f_{i,(-1)}\partial^{j_i}\xi_i$.  The normal ordering formula
 (\ref{norm-order-form-module}) gives
 \begin{eqnarray}
 & &g_{(-1)}f_{i,(-1)}\partial^{j_i}\xi_i=\nonumber\\
 & &(gf_i)_{(-1)}\xi_i-\sum_{j=0}^\infty f_{i,(-j-2)}g_{(j)}\xi_i-\sum_{j=0}^\infty  g_{(-j-2)}f_{i,(j)}\xi_i=\nonumber\\
 & &(gf_i)_{(-1)}\xi_i-\sum_{j=0}^\infty\frac{1}{j!} \partial^j f_{i,(-1)}g_{(j)}\xi_i-\sum_{j=0}^\infty  \frac{1}{j!} \partial^j g_{(-1)}f_{i,(j)}\xi_i,
 \nonumber
 \end{eqnarray}
 which is determined by the considerations above.

 Therefore,  the room for maneuver is only provided by the Lie* bracket on $\CL^{ch}_A$.  If $[.,.]$ is one such bracket, then
any bracket is
\[
[\xi.,.\tau]_\alpha=[\xi,\tau]+\alpha(\xi,\tau) \text{ for some } \alpha\in P^*_{\{1,2\}}(\{J_\infty T_A,J_\infty T_A\}, J_\infty A).
\]

It easily follows from (\ref{what-l-linear-means}) that $\alpha$ must be $J_\infty A$-linear,  see (\ref{defi-multil-in shrie-ca}) for the definition of  $J_\infty A$-linearity. The antisymmetry of a Lie* bracket implies that $\alpha$ must be antisymmetric.
 The Jacobi identity,
 \[
 [\xi_1[\xi_2,\xi_3]_\alpha]_\alpha- [\xi_2[\xi_1,\xi_3]_\alpha]_\alpha- [[\xi_1,\xi_2]_\alpha,\xi_3]_\alpha=0,
 \]
 implies
\begin{multline}
\nonumber
[\xi_1,\alpha(\xi_2,\xi_3)]-[\xi_2,\alpha(\xi_1,\xi_3)]+[\xi_3,\alpha(\xi_1,\xi_2)]-\\
\alpha(\sigma[\xi_1,\xi_2],\xi_3)+\alpha(\sigma[\xi_1,\xi_3],\xi_2)-\alpha(\sigma[\xi_2,\xi_3],\xi_1)=0
\end{multline}
Since the Lie* bracket $[.,.]$ restricted to $J_\infty A$ is the pull-back of the 
canonical action of $J_\infty T_A$ on $J_\infty A$ (item (ii) of the definition in sect.~\ref{sh-ex-seq-def-chir-alg}),
this means that $\alpha$ is a closed 2-cochain of $J_\infty T_A$ with coefficients in $J_\infty A$, which satisfies an extra condition of being $J_\infty A$-linear,
see the definition of the  Chevalley complex in sect.~\ref{lie-star-chevalley}.  

More generally, define $C^n_{J_\infty A}(J_\infty T_A, J_\infty A)\subset C^n(J_\infty T_A, J_\infty A)$ to be the subspace of $J_\infty A$-linear operations. It is easy to see that $C^\bullet_{J_\infty A}(J_\infty T_A, J_\infty A)\subset C^\bullet(J_\infty T_A, J_\infty A)$ is a subcomplex, and
as such it is called the {\em Chevalley--De Rham complex}; this definition makes sense for any Lie* algebroid, \cite{BD}, 1.4.14.

To conclude, given a chiral algebroid $\CL^{ch}_A$ and $\alpha\in C^{2,cl}_{J_\infty A}(J_\infty T_A, J_\infty A)$ we have defined another chiral algebroid, to be denoted
$\CL^{ch}_A(\alpha)$; furthermore, any chiral algebroid is isomorphic to $\CL^{ch}_A(\alpha)$ for some $\alpha$.

The description of morphisms is also similar  to sect.~\ref{defn-p-lie}.  By definition, each morphism must have the form 
\[
\xi\mapsto \xi +\beta(\xi),\; \beta\in 
C^1_{J_\infty A}(J_\infty T_A, J_\infty A)=\text{Hom}_{J_\infty A}(J_\infty T_A, J_\infty A).
\]
A quick computation, no different from the ordinary case,  will show
\[
\text{Hom}_{Ch-Alg}(\CL^{ch}_A(\alpha_1),\CL^{ch}_A(\alpha_2))=\{\beta\in C^1_{J_\infty A}(J_\infty T_A, J_\infty A)\text{ s.t. }d\beta=\alpha_1-\alpha_2\}.
\]
\begin{sloppypar}
This can be rephrased as follows. Let $C^{[1,2>}(J_\infty T_A)$ be the category with objects $C^{2,cl}_{J_\infty A}(J_\infty T_A, J_\infty A)$ and morphisms
$\text{Hom}(\alpha_1,\alpha_2)=\{\beta\in C^1_{J_\infty A}(J_\infty T_A, J_\infty A)\text{ s.t. }d\beta=\alpha_1-\alpha_2\}$. 
If $T_A$ is a free $A$-module, then the category of chiral $A$-algebroids is a $C^{[1,2>}(J_\infty T_A)$-torsor. It is non-empty if $T_A$ has a finite abelian basis; this follows from sect.~\ref{ex-poly-ring}.
\end{sloppypar}

\subsection{ }
\label{localaz-chir-alg}
These considerations can be localized in an obvious manner. For any smooth $X$, one obtains a tangent Lie* algebroid $T^{ch}_X$ and 
 a gerbe of chiral algebroids over $J_\infty X$, bound by the 
complex $C^1_{J_\infty X}(T^{ch}_X,\CO_{J_\infty X})\rightarrow C^{2,cl}_{J_\infty X}(T^{ch}_X,\CO_{J_\infty X})$. This gerbe is locally non-empty, as follows
from sect.~\ref{ex-poly-ring}. The calculation of its characteristic class, in this and much greater generality, can be found in \cite{BD}, 3.9.22. We shall review below (sect.~\ref{gorb-mal-schechtm})
the case of a graded chiral agebroid.

\subsection{ }
\label{univ-envel-chir-alg}
The chiral enveloping algebra $U^{ch}(\CL_A^{ch})$ attached to $\CL_A^{ch}$ if the latter is regarded as a Lie* algebra,  sect.~\ref{vert-envel-lie-star}, does not ``know'' about the
chiral structure that $\CL_A$ carries. This leads to the existence of a canonical ideal as follows. Consider two elements
$a_{(-n)}\xi\in\CL_A^{ch}$ and $a_{[-n]}\xi\in\CL_A\cdot U^{ch}(\CL_A^{ch})$, $n>0$, where $a\in J_\infty A\subset U^{ch}(\CL_A^{ch})$ and $\xi\in\CL_A^{ch}\subset U^{ch}(\CL_A^{ch})$.  Since both these products, $_{(-n)}$, reflecting the chiral $J_\infty A$-module structure of
$\CL_A^{ch}$, and $_{[n]}$, reflecting the chiral algebra structure of $U^{ch}(\CL_A^{ch})$, sect.~\ref{vert-envel-lie-star}, satisfy the same Borcherds commutator formula, cf.
(\ref{borch-in-hl}) and (\ref{what-l-linear-means}), their difference satisfies
\[
\tau_{[m]}(a_{(-n)}\xi-a_{[-n]})\xi\in\CL_A^{ch}\cdot U^{ch}(\CL_A^{ch}))=0\text{ if } m\geq 0.
\]
This is a familiar {\em singular vector} condition. Denote by  $J$ the maximal chiral ideal of $U^{ch}(\CL_A^{ch})$ generated by all such elements along with the difference
$1_{A}-1_{U^{ch}}$, where $1_A\in  J_\infty A$ and $1_{U^{ch}}\in U^{ch}(\CL_A^{ch})$ are units.  It is practically obvious that $ D^{ch}_{\CL_A^{ch}}$ defined to be the quotient
$U^{ch}(\CL_A^{ch})/J$ is a CDO over $A$, sect.~\ref {defftn-cdo}, at least if the tangent algebroid  $T_A$ is  a (locally) free $A$-module. In fact, the assignment
$\CL_A\mapsto D^{ch}_{\CL_A^{ch}}$ is an equivalence of categories.

All of this is, of course,  parallel to sect.~\ref{lie-alg-envelope}.

{\em Example.} If $\CL_A^{ch}$ is $F_1 A[x_{ij},\partial_{mn};\;j,n<0]$ introduced in sect.~\ref{ex-poly-ring}, then $ D^{ch}_{\CL_A^{ch}}=A[x_{ij},\partial_{mn};\;j,n<0]\}$.

\subsection{ }
We shall say that  a chiral algebra $\CV$ is $\BZ$-{\em graded} if $\CV=\oplus_{n\in\BZ}\CV_n$ s.t.
$\CV_{n(j)}\CV_m\subset \CV_{m+n-j-1}$ and  $\partial(\CV_n)\subset\CV_{n+1}$.  A similar definition also applies to coisson algebras, sect.~\ref{comm-modul-act-prod-comm}. Here is the origin of this concept.

Let $L$ be a Lie* algebra. We say that $L$ acts on a chiral algebra $\CV$ if $\CV$ is an $L$-module such that 
the chiral bracket $\mu\in P^{ch}_{\{1,2\}}(\{\CV,\CV\},\CV)$ is $L$-linear, cf. sect.~\ref{defn-chir-algebr} (v) and Remark (1).

Let $\CV ec$ be a free $R=\BC[\partial]$-module on 1 generator $l$. Make it into a Lie* algebra  by defining a Lie* bracket so that
\[
l\otimes l\mapsto -l\partial\otimes 1+2l\otimes \partial_1.
\]
This is equivalent to saying that $l_{(0)}l=-l\partial$, $l_{(1)}l=2l$.  Call an action of $\CV ec$ on $\CV$ {\em nice} if $l^{\CV}_{(-1)}v=-v\partial$. For example, the adjoint
action of $L$ on itself is nice. 

One readily verifies that  $\CV$ is  $\BZ$-graded iff $\CV$ carries a nice action of $\CV ec$ such that the operator $l_{(1)}^{\CV}$ is diagonalizable.  The equivalence is established by stipulating $\CV_n=\{v\text{ s.t. } l^{\CV}_{(1)}v=nv\}$.

Notice the (easy to verify) isomorphism $h(\CV ec[x])\iso T_{\BC[x]}$ defined by $l\otimes x^n\mapsto -x^n \partial/\partial x$, sect.~\ref{Lie-alf-fourier-comp}, and so the grading
operator has the meaning of $-x\partial/\partial x$.

\subsection{ }
\label{gorb-mal-schechtm}
Now it should be clear what a $\BZ$-graded chiral algebroid is; we call it $\BZ_+$-graded if $(\CL^{ch}_A)_n=\{0\}$ provided $n<0$. Classification
of  $\BZ_+$-graded chiral algebroids is simpler and more explicit, \cite{GMS}. We continue under the assumption that $T_A$ is a free $A$-module
with a finite {\em abelian } basis $\{\tau_i\}$. Denote by $\{\omega_i\}\subset\Omega_A$ the dual basis: $\omega_i(\tau_j)=\delta_{ij}$. In this case there is
always an $\CL^{ch}_A$ determined by the requirements $\tau_{i(n)}\tau_j=0$ if $n\geq 0$, see sect.~\ref{ex-poly-ring}.

Notice that the quasiclassical object $J_\infty A\oplus J_\infty T_A$ is naturally $\BZ_+$-graded: place $A\subset J_\infty A$ in degree 0, $T_A\subset J_\infty T_A$ in degree 1, 
and use the fact that
$\partial$ has degree 1.  We seek, therefore, a classification of those  $\BZ_+$-graded chiral algebroids whose grading induces the indicated one on the quasiclassical object.

Having split $\CL^{ch}_A$ into the direct sum $J_\infty A\oplus J_\infty T_A$ as in sect.~\ref{class-chir-algebr}, we obtain that a variation of the Lie* bracket 
is an operation $\alpha(.,.)\in P^*_{\{1,2\}}(\{J_\infty T_A,J_\infty T_A\},J_\infty A)$, which is $J_\infty A$-bilinear. Since $J_\infty T_A$ is a free $\BC[\partial]$ module, it is determined
by its values on $T_A\subset J_\infty T_A$:
\[
\alpha(\xi,\eta)=\sum_{n=0}^\infty \alpha_n(\xi,\eta)\otimes \frac{\partial_1^n}{n!},\; \xi,\eta\in T_A.
\]
The grading condition demands that at most   2 components may be nonzero:
\[
\alpha(\xi,\eta)= \alpha_0(\xi,\eta)\otimes 1+ \alpha_1(\xi,\eta)\otimes \partial_1,
\]
where $\alpha_0(\xi,\eta)\in\Omega_A= (J_\infty A)_1$, $\alpha_1(\xi,\eta)\in A=( J_\infty A)_0$.  Furthermore, varying the splitting $J_\infty T_A\hookrightarrow (\CL_A^{ch})_1$
by sending $\tau_i\mapsto\tau_i-1/2\sum_j\alpha_1(\tau_i,\tau_j)\omega_j$ ensures that $\alpha_1$ is 0.

Component $\alpha_0$, as it stands,  is an antisymmetric $A$-bilinear map from $T_A$ to $\Omega_A$, hence  $\alpha_0\in\Omega^2_A\otimes_A\Omega_A$. 
The relation $\xi_{(0)}(\eta_{(1)}\gamma)=(\xi_{(0)}\eta)_{(1)}\gamma+\eta_{(1)}(\xi_{(0)}\gamma)$,  which is (\ref{jacobi})  with $n=0$, $m=1$, shows that in fact
$\alpha_0$ is totally antisymmetric and so belongs to $\Omega^3_A$. Finally, the relation $\xi_{(0)}(\eta_{(0)}\gamma)=(\xi_{(0)}\eta)_{(0)}\gamma+\eta_{(0)}(\xi_{(0)}\gamma)$,  which is (\ref{jacobi})  with $n=m=0$,  shows that $\alpha_0$ is, moreover, a closed 3-form.

Similarly, a change of splitting $\xi\mapsto\xi+\beta(\xi)$ preserves the grading precisely when $\beta\in\Omega_A\otimes_A\Omega_A$ and the normalization we
chose ($\tau_{i(1)}\tau_j=0$) requires that $\beta\in\Omega^2_A$. The effect of this on $\alpha_1$ is $\alpha_1\mapsto \alpha_1-d_{DR}\beta$.

More formally, the meaning of these computations is as follows.  The truncated Chevalley-- De Rham complex  $C^1_{J_\infty A}(J_\infty T_A,J_\infty A)\longrightarrow C^{2,cl}_{J_\infty A}(J_\infty T_A,J_\infty A)$, introduced in sect.~\ref{class-chir-algebr}, is graded, and its degree 0 component,
$C^1_{J_\infty A}(J_\infty T_A,J_\infty A)[0]\longrightarrow C^{2,cl}_{J_\infty A}(J_\infty T_A,J_\infty A)[0]$, describes the category of $\BZ_+$-graded chiral algebroids. 
Described above is a map of the truncated De Rham complex $\Omega^2_A\longrightarrow\Omega^{3,cl}_A$ to $C^1_{J_\infty A}(J_\infty T_A,J_\infty A)[0]\longrightarrow C^{2,cl}_{J_\infty A}(J_\infty T_A,J_\infty A)[0]$; e.g., this map sends 
\begin{eqnarray}
\nonumber
\Omega^{3,cl}_A\ni\alpha_0&\mapsto&  \alpha \in C^{2,cl}_{J_\infty A}(J_\infty T_A,J_\infty A)[0]\text{ s.t. }\alpha(\xi,\eta)= \alpha_0(\xi,\eta,.)\otimes 1,\\
\Omega^{2}_A\ni\beta_0&\mapsto&  \beta \in C^{1}_{J_\infty A}(J_\infty T_A,J_\infty A)[0]\text{ s.t. }\alpha(\xi)= \alpha_0(\xi,.)\otimes 1.\nonumber
\end{eqnarray}
Analogously to $C^{[1,2>}(J_\infty T_A)$, sect.~\ref{class-chir-algebr}, introduce $\Omega^{[2,3>}_A$, the category with objects $\Omega^{3,cl}_A$ and morphisms
$\text{Hom}(\alpha_1,\alpha_2)=\{\beta\in\Omega^2_A\text{ s.t. }\alpha_1-\alpha_2=d_{DR}\beta\}$.  The map of complexes just defined gives a functor
$\Omega^{[2,3>}_A\longrightarrow C^{[1,2>}(J_\infty T_A)[0]$.  The point is:  {\em  this functor is an equivalence of categories.}

To summarize:  {\em if $A$ is such that $T_A$ is a free $A$-module with a finite abelian basis, then the category of chiral $A$-algebroids is a $\Omega^{[2,3>}_A$-torsor.}

\subsection{ }
\label{ch2}
These considerations can be localized so as to obtain, over any smooth $X$, a gerbe of $\BZ_+$-graded CDOs bound by the complex 
$\Omega^2_X\longrightarrow\Omega^{3,cl}_X$; this gerbe is locally non-empty. Its characteristic class is $ch_2(T_X)$. The details of this computation can be found in \cite{GMS}; cf. \cite{BD}, 3.9.23.

\subsection{ }
\label{somewhat-tw-cdo}
One can slightly relax the $\BZ_+$-graded condition by demanding that the CDO be filtered, i.e., that
\[
[J_\infty T_A,J_\infty A]\subset J_\infty T_A\otimes 1\oplus \Omega_A\otimes 1\oplus A\otimes 1\oplus A\otimes \partial_1,
\]
here the summand $A\otimes 1$ is the one that was prohibited in sect.~\ref{gorb-mal-schechtm}. In other wards, we allow variations of the form
\[
[\xi,\eta]_{\alpha,\beta}=[\xi,\eta]+\alpha(\xi,\eta)+\beta(\xi,\eta),\; \alpha(\xi,\eta)\in\Omega_A,\, \beta(\xi,\eta)\in A.
\]
Just as before, one obtains $\alpha(.,.)\in\Omega^{3,cl}_A$, $\beta(.,.)\in \Omega^{2,cl}$, and (provided $T_A$ has an abelian basis) the category
of filtered CDOs is an $\Omega^{[2,3>}_A\times \Omega^{[1,2>}_A$-torsor, thereby getting a cross between the Picard-Lie (sect.~\ref{defn-p-lie}) and graded chiral
algebroid. This is similar to but different from the concept of a {\em twisted} CDO introduced (and used) in \cite{AChM,AM}. On the other hand,
examples of such CDOs have already crept in the literature: \cite{H,LinMath}

\bigskip

\section{chiral $\infty$-algebroids}
\label{chiral-infty-algebro}
 The main result is
Theorem~\ref{class-oo-chi-alge}, which is very similar to Lemma~\ref{res-class-picardlie-infty} except that the ordinary (derived) De Rham complex is
replaced with its version in the world of the Beilinson-Drinfeld pseudo-tensor category.

\subsection{ }
\label{descr-obf-chiral-infty-al}

A chiral $\infty$-algebroid over a DGA $A$ is a short exact sequence, cf. sect.~\ref{defn-chir-algebr},
\begin{equation}
\label{ex-seq-ari-inft-chiral-alge}
0\longrightarrow J_\infty A\stackrel{\iota}{\longrightarrow}\CL_A^{ch}\stackrel{\sigma}{\longrightarrow} J_\infty T_A\longrightarrow 0,
\end{equation}
where 

$A$ is a commutative finitely generated DGA with degree 1 differential $D$;

 $J_\infty A$ and $J_\infty T_A$
 are as in {\em loc. cit.}, except that they carry an extra differentia $D$; in particular,  $J_\infty A$
is a commutative  DG chiral algebra and $J_\infty T_A$  is a DG  Lie* $J_\infty A$-algebroid (sect.~\ref{comm-modul-act-prod-comm});

 $\CL^{ch}_A$ is a DG $\text{Lie}^*_\infty$ algebra with operations $l_n\in P^*_{[n]}(\{\CL^{ch}_A\},\CL^{ch}_A)$, $\text{deg}l_n=2-n$, $n\geq 1$, see sect.~\ref{lie-star-infty},
 and a chiral DG $J_\infty A$-module, which is  defined by an operation $\mu\in P^{ch}_{[2]}(\{J_\infty A, \CL^{ch}_A\},\CL^{ch}_A)$.

The following conditions must hold:

(i) The morphism $ J_\infty A\stackrel{\iota}{\longrightarrow}\CL_A^{ch}$ is a  DG chiral $J_\infty A$-module and a strict  $\text{Lie}_\infty^*$ algebra morphism.

(ii) If we let $\mu^*\in P^{*}_{[2]}(\{J_\infty A, \CL^{ch}_A\},\CL^{ch}_A)$ be the operation determined by $\mu$ via (\ref{how-chira-ope-def-star-opera}), then $\mu^*=l_2|_{J_\infty A\otimes\CL^{ch}_A}$.

(iii) The morphism $ \CL_A^{ch}\stackrel{\sigma}{\longrightarrow} J_\infty T_A$  is a  DG chiral $J_\infty A$-module and a strict  $\text{Lie}^*_\infty$ algebra morphism.

(iv) By (i) and (iii) $J_\infty T_A=\CL^{ch}_A/J_\infty A$ operates on $J_\infty A$ as a Lie* algebra. We require that this action coincide with the tautological action
of  $J_\infty T_A$ on $J_\infty A$.

(v) The operation $l_2$ is a derivation of the chiral action $\mu$, cf (\ref{meani-of-lch-linea-defni-cdo}).Namely
\[
l_2(.,\mu(.,.))=\mu(l_2(.,.),.)+\mu(.,l_2(.,.))^{(1,2)}.
\]

(vi) The operations $l_n$, $n\geq 3$ are $J_\infty A$-valued and factor through  the morphism $ \CL_A^{ch}\stackrel{\sigma}{\longrightarrow} J_\infty T_A$.
The corresponding operations, to be also denoted $l_n\in P^*_{[n]}(\{J_\infty T_A\},J_\infty A)$, are $J_\infty A$-multilinear, as defined in (\ref{defi-multil-in shrie-ca}).

\subsubsection{ } Of course, an ordinary chiral algebroid with differential is an example of a chiral $\infty$-algebroid. The Feigin-Semikhatov construction
discussed in the introduction gives us an example, where $A=\BC[x,\xi]$, $x$ even, $\xi$ odd, 
$\CL^{ch}_A=J_\infty\BC[x,\xi]\oplus J_\infty T_{\BC[x,\xi]}$ and
$D= (x^m\partial_\xi)_{(0)}$, cf. sect.~\ref{ex-poly-ring}.

\subsection{ }
\label{chira-inft-alge-morphis}

Define a morphism of chiral  $\infty$-algebroids, $ (\CL^{ch}_A,\iota,\sigma)\rightarrow(\tilde{\CL}_A^{ch},\tilde{\iota},\tilde{\sigma})$ to be a morphism of $\text{Lie}^*_\infty$ algebras
$f=\{f_n\}: \CL^{ch}_A\rightarrow \tilde{\CL}_A^{ch}$ (as defined in sect.~\ref{defni-morph-listarinfty-case}) that satisfies the following 2 conditions:

(i) each $f_n\in SP^*_{[n]}(\{\CL^{ch}_A[1]\}, \tilde{\CL}^{ch}_A[1])$, is $J_\infty A$ $n$-linear; furthermore, if $n>1$, then $f_n$ 
is $J_\infty A$-valued and factors through the 
 map
$\CL_A^{ch}\stackrel{\sigma}{\longrightarrow} J_\infty T_A$, i.e., vanishes if one of the arguments is in $J_\infty A$;

(ii) the component  $f_1\in P^*_{[1]}(\{\Delta^{(\CS)}\CL^{ch}_A[1]\}, \Delta^{(\CS)}\tilde{\CL}^{ch}_A[1])$, which according to (\ref{exte-morphi-2})
can be regarded as a morphism $f_1:\CL^{ch}_A\rightarrow\tilde{\CL}_A^{ch}$ (cf. (\ref{expl-descr-morph-infty}),
makes the following a commutative diagram of chiral $J_\infty A$-module morphisms:

\begin{displaymath}
\label{comm-diagr-morph-chi-alge-inf}
\xymatrix{ 0\ar[r]&A\ar[r]^{\iota}&\CL_A\ar[r]^{\sigma}\ar[d]^{f_1}&
J_\infty T_A\ar[r]&0\\
 0\ar[r]&A\ar@{=}[u] \ar[r]^{\tilde{\iota}}& \tilde{\CL}_A^{ch}\ar[r]^{\tilde{\sigma}}&J_\infty T_A\ar@{=}[u]\ar[r]&0}
 \end{displaymath}
 
 Note that according to (i)   $f_n$ can be regarded as an element of $SP^*_{[n]}(\{J_\infty T_A[1]\}, J_\infty A[1])$ if $n>1$.
 
 \subsection{ }
 \label{cate-of-forms-bd-settti}
 Recalled in sect.~\ref{class-chir-algebr}, the De Rham-Chevalley complex $(C^\bullet_{J_\infty A}(J_\infty T_A,J_\infty A),d)$ in the present situation acquires
 an extra differential, induced by $D$, the differential of $A$, and an extra grading, also inherited from $A$:
 $C^\bullet_{J_\infty A}(J_\infty T_A,J_\infty A)=\oplus_{n\leq 0}C^\bullet_{J_\infty A}(J_\infty T_A,J_\infty A)[n]$.
  Denote by $LC^\bullet_{J_\infty A}(J_\infty T_A,J_\infty A)$ the completion in the De Rham direction of the corresponding total complex; one has
  $LC^n_{J_\infty A}(J_\infty T_A,J_\infty A)=\prod_k C^k_{J_\infty A}(J_\infty T_A,J_\infty A)[n-k]$, the differential being $d+D$; this is a straightforward analogue of the
  Illusie construction \cite{Ill1,Ill2}.
  
  Denote by $LC^{[1,2>}_{J_\infty A}(J_\infty T_A,J_\infty A)$ the category with objects $LC^{2,cl}_{J_\infty A}(J_\infty T_A,J_\infty A)$ and 
  \[
  Hom(\beta,\gamma)=\{\alpha\in LC^{1}_{J_\infty A}(J_\infty T_A,J_\infty A)\text{ s.t. }\beta-\gamma=(d+D)\alpha\},
  \]
   cf. sect.~\ref{inftya-alg-categ}.
 
 \begin{thm}
 \label{class-oo-chi-alge}
 Let $(A,D)$ be a finitely generated polynomial commutative DGA, the degree of $D$ being 1.
 The category of chiral $\infty$-algebroids over $A$ is a torsor over $LC^{[1,2>}_{J_\infty A}(J_\infty T_A,J_\infty A)$.
 \end{thm}

{\em Proof.}
Given  a $(\CL^{ch}_A,\{l_n\})$ , any other chiral $\infty$-algebroid can be defined by varying $l_n\mapsto l_n+\alpha_n$,
for some $\alpha_n\in  P^*_{[n]}(\{J_\infty T_A,J_\infty A)$, and by definition
\[
\alpha_n\in C^n_{J_\infty A}(J_\infty T_A,J_\infty A)[2-n].
\]
The quadratic relations that appear in the definition of a $\text{Lie}^*_\infty$ algebra, sect.~\ref{lie-star-infty}, are equivalent to the cocycle condition
 \[
 \{\alpha_n\}\in LC^{2,cl}_{J_\infty A}(J_\infty T_A,J_\infty A);
 \]
this discussion is to to sect.~\ref{class-chir-algebr} exactly what sect.~\ref{inftya-alg-categ}
is to sect.~\ref{defn-p-lie}. In fact, the derivation of the above cocycle condition from the $\text{Lie}^*_\infty$ algebra definition is no different
from the corresponding proof in sect.~\ref{inftya-alg-categ}. This defines an action of $LC^{[1,2>}_{J_\infty A}(J_\infty T_A,J_\infty A)$ on the category of
chiral $\infty$-algebroids, $(\CL^{ch}_A,\{l_n\}),\{\alpha_n\}\mapsto (\CL^{ch}_A,\{l_n+\alpha_n\})$.

\begin{sloppypar}
By definition, see sect.~\ref{chira-inft-alge-morphis}, morphisms are  collections of degree 0 operations $\hat{\beta}_n: SP^*_{[n]}(J_\infty T_A[1], J_\infty  A[1])$. 
The actual morphism that such a collection defines operates as follows, cf. (\ref{expl-descr-morph-infty-in-bd-ca}):
\end{sloppypar}
\[
f(x)=x+\hat{\beta}_1(x),
\]
\[
f(x_1,x_2)=(x_1+\hat{\beta}_1(x_1))(x_2+\hat{\beta}_1(x_2))+\hat{\beta}_2(x_1,x_2),
\]
etc.  Such morphisms are automatically automorphisms; in fact, $f^{-1}$ is the morphism defined by the collection $\{-\hat{\beta}_n\}$.

The effect $f$ has on the coderivation $\hat{l}=\sum_n\hat{l}_n$ is this: $\hat{l}\mapsto f\circ\hat{l}\circ f^{-1}$.
To compute the difference, $\hat{l}-f\circ\hat{l}\circ f^{-1}$, remove the hats by defining 
$\beta_n=s\circ\hat{\beta}_n\circ((s^{-1})^{\otimes n})$, cf. sect.~\ref{using-RS-symm-alg-to-descri-star-infty-li}. We have $\beta_n\in C_{J_\infty A}(J_\infty T_A,J_\infty A)^n_A[1-n]$, as 
in {\em loc. cit.}.
A  straightforward computation will then reveal that 
\[
\hat{l}-f\circ\hat{l}\circ f^{-1}= (d_{DR}+\text{Lie}_\partial)(\{\beta_n\})\hat{ },
\]
as desired. $\qed$


\bigskip


\end{document}